\documentclass[a4paper]{article}

\usepackage{amsmath}
\usepackage{latexsym}
\usepackage{graphicx}
\usepackage{amsfonts}
\makeatletter
\def\@maketitle{\newpage
    \null
    \vskip .8truein
    \begin{center}%
     {\bf \@title \par}%
     \vskip 1.5em
     {\small
      \lineskip .5em
      \begin{tabular}[t]{c}\@author
      \end{tabular}\par}%
    \end{center}%
    \par
    \vskip .4truein}
\@addtoreset{equation}{section} \@addtoreset{theorem}{section}
\@addtoreset{lemma}{section} \@addtoreset{proposition}{section}
\@addtoreset{definition}{section} \@addtoreset{corollary}{section}
\@addtoreset{remark}{section}

\let\O=\Omega

\let\d=\delta

\let\s=\sigma

\let\nn=\nonumber

\let\de=\partial 
\let\ol=\overline
\newcommand{\re}{{{I\!\!R}}}
\let\ul=\underline

\def\R{\mathbb{R}}
\let\a=\alpha
\let\d=\delta
\let\e=\varepsilon
\let\g=\gamma
\let\l=\lambda
\let\vf=\varphi

\def\mX{\mathcal X}

\newtheorem{theorem}{Theorem}[section]
\newtheorem{lemma}{Lemma}[section]
\newtheorem{proposition}{Proposition}[section]
\newtheorem{definition}{Definition}[section]
\newtheorem{corollary}{Corollary}[section]
\newtheorem{remark}{Remark}[section]
\newtheorem{example}{Example}[section]
\DeclareMathOperator{\argmax}{argmax}

\def\proof{\list{}{\setlength{\leftmargin}{0pt}
                      \parskip=0pt\parsep=0pt\listparindent=2em
                      \itemindent=0pt}\item[]\futurelet\testchar\@maybe}

\def\@maybe{\ifx[\testchar \let\next\@Opt
          \else \let\next\@NoOpt \fi \next}
\def\@Opt[#1]{{\it Proof of #1.\ }}\def\@NoOpt{{\it Proof.\ }}

\title{\Large \bf Comparison principles and Dirichlet problem for equations of Monge-Amp\`ere type associated to vector fields
\footnotemark[1]
}

\author{{\large \sc Martino Bardi, Paola Mannucci}\\
 \rm Dipartimento di Matematica Pura e Applicata,\\
Universit\`a degli Studi di Padova,\\
Via Trieste, 63, 35121, Padova, Italy,\\
bardi@math.unipd.it, mannucci@math.unipd.it
 }

\begin{document}

\nonstopmode

\maketitle
\renewcommand{\thefootnote}{\fnsymbol{footnote}}
\footnotetext[1]{Work partially supported by
the Italian M.I.U.R. project "Viscosity,
metric, and control theoretic methods for nonlinear partial
differential equations''.}



\begin{abstract} We study 
partial differential equations of Monge-Amp\`ere type involving the derivates with respect to a family $\mX$ of vector fields of Carnot type. The main result is a comparison principle among viscosity subsolutions, convex with respect to $\mX$, and viscosity supersolutions (in a weaker sense than usual), which implies the uniqueness of solution to the Dirichlet problem. Its assumptions include the equation of prescribed horizontal Gauss curvature in Carnot groups. By Perron method we also prove the existence of a solution either under a growth condition of the nonlinearity with respect to the gradient of the solution, or assuming the existence of a subsolution attaining continuously the boundary data, therefore generalizing some classical result for Euclidean Monge-Amp\`ere equations.
\end{abstract}



\tableofcontents

\section{Introduction}

For a given family of $C^{1,1}$ vector fields $\mX=\{X_1, ..., X_m\}$ in $\R^n$, $m\leq n$, the $\mX$-gradient and  symmetrized $\mX$-Hessian matrix of a function $u$ are 
\[
D_{\mathcal X}u := (X_1u, ..., X_mu),
\qquad
(D^2_{\mathcal X}u)_{ij} :=\left( X_i X_j u + X_j X_i u\right) /2 .
\]
The main examples we have in mind are the vector fields that generate a homogeneous Carnot group \cite{Bell, BLU}, and in that case $D_{\mathcal X}u$  and 
$D^2_{\mathcal X}u$ are called, respectively, the horizontal gradient and  the horizontal Hessian.
We consider fully nonlinear partial differential equations of the form
\begin{equation}
\label{MAh}
- \det 
D_ {\mathcal X}^2 u
 + H(x, u, D_{\mathcal X}u)=0, \ \mbox{in}\ \Omega,
\end{equation}
where $\Omega\subseteq \re^n$ is open 
and bounded
and $H$ is 
at least continuous and nonnegative.  
In  the case when the vector fields are the canonical basis of $\R^n$, that we call the Euclidean case, this is a classical equation of Monge-Amp\`ere type. 
 We recall that in the Euclidean case  the Monge-Amp\`ere equations are elliptic on convex functions.
These equations arise in several 
 problems, mostly of differential geometry, and have a wide literature, especially on the regularity of solutions: see, e.g.,  the books \cite{Po78, GT, Ba94, Au, Kry, Gbook} and the papers \cite{CY77, CY82, Del, CNS, PLL, PLL2, Kry98, Tru90, Tru06}. For the recent applications to optimal transportation problems we refer to \cite{ACBBV, Vbook1, Ca}  and the references therein.
 

Partial differential equations with an elliptic structure relative to vector fields that do not span the whole space $\R^n$ are often called subelliptic, see, e.g., the recent book of Bonfiglioli, Lanconelli, and Uguzzoni \cite{BLU} for a comprehensive survey of the linear theory.
A theory of fully nonlinear subelliptic equations was started a few years ago by Bieske \cite{B1, B2} and Manfredi   \cite{Man, BBM} using viscosity methods, and the Monge-Amp\`ere equation 
\begin{equation}
\label{MAf}
- \det 
D_ {\mathcal X}^2 u
 + f(x) = 0, \ \mbox{in}\ \Omega,
\end{equation}
was listed among the main examples, with $X_1,...,X_m$ generators of a given Carnot group.
Moreover, a number of authors studied in the last five years several notions of convexity in Carnot groups \cite{LMS, DGN, GutMon, GutMon2, BaRic, DGNT, W2, Mag, Ric, GarTou, JLMS}, and one of their motivations was the connection with Monge-Amp\`ere equations on such groups. However, 
little is known 
about them so far. 
We mention the Comparison Principle among smooth sub- and supersolutions of   \eqref{MAf} in the Heisenberg group proved by Gutierrez and Montanari \cite{GutMon} (among other results).

The dependence on the gradient $D_{\mathcal X}u$ in $H$ is motivated by various possible applications.
A first interesting example is the subelliptic analogue of the prescribed Gauss curvature equation. In fact, Danielli, Garofalo and Nhieu \cite{DGN} defined the horizontal Gauss curvature of the graph of the smooth function $u$ on a Carnot group as
\[
K_h(x) := \det (D_{\mathcal X}^2 u) \left(1+|D_{\mathcal X} u|^2\right)^{-\frac{m+2}{2}},
\]
so \eqref{MAh} becomes the prescribed horizontal Gauss curvature equation if 
\begin{equation}
\label{Gauss}
H(x, r, q) = k(x) \left(1+|q|^2\right)^{\frac{m+2}{2}}
\end{equation}
for a given continuous $k \,:\, \ol\Omega \to [0, +\infty[$. A different hypoelliptic Monge-Amp\`ere-type equation was proposed in \cite{Sto} for a financial problem. Finally, the extension of the theory of optimal transportation to the realm of Sub-Riemannian manifolds was started recently by \cite{AR} and \cite{FiRi} and it might lead to equations of the form  \eqref{MAh}, or variants of it.

This paper is devoted to a study of the subelliptic Monge-Amp\`ere-type equations \eqref{MAh} within the theory of viscosity solutions, see \cite{CIL, CC, CCime, BCD}. In particular we establish the well-posedness of the Dirichlet problem under rather general conditions. The first part of the paper deals with comparison results among sub- and supersolutions, and the second part with the existence of solution by the Perron method.

 The new difficulties we encounter for the Comparison Principles are three.

1.  The PDE \eqref{MAh} is degenerate elliptic only on %
convex  functions with respect to the vector fields $X_1, ..., X_m$,  briefly {\em ${\mathcal X}$-convex}. We say that an u.s.c. function $u$ on $\ol\Omega$ is ${\mathcal X}$-convex  if it satisfies
$
-D^2_{\mathcal X}u\leq 0
$
in $\Omega$ in viscosity sense, that is, 
\begin{equation}
\label{hconv-test}
D^2_{\mathcal X}\vf(x)\geq 0 \quad \forall\; \vf\in C^2(\Omega), x\in\argmax(u-\vf).
\end{equation}
This notion was intruduced by  Lu, Manfredi, and Stroffolini \cite{LMS} for the Heisenberg group under the name of  v-convexity. It was  extended very recently to general $C^2$ vector fields by the first author and Dragoni \cite{BD}, who proved the equivalence with the convexity along trajectories of the fields.
In the case of Carnot groups it coincides with the 
 geometric notion of {\em horizontal convexity}, a fact proved under different assumptions by several authors  \cite{LMS, BaRic, W2, Mag, Ric, JLMS}, see also \cite{DGN}  for connections with other notions.
Our comparison results will concern a ${\mathcal X}$-convex viscosity subsolution of \eqref{MAh} and a viscosity supersolution defined with strictly ${\mathcal X}$-convex test functions. This is inspired by the treatment of the Euclidean case by Ishii and P.-L. Lions  \cite{IL} and is equivalent to 
comparing 
  sub- and supersolutions of 
\[
\max\{-\lambda_{min}(D_ {\mathcal X}^2 u), - \det D_ {\mathcal X}^2 u+ H(x, u, D_{\mathcal X}u)\}=0, \ \mbox{in}\ \Omega,
\]
where $\lambda_{min}$ denotes the minimal eigenvalue. 

In the classical case convex functions are locally Lipschitz continuous with respect to the Euclidean norm, so 
there is an interior  gradient bound for the subsolution.
The corresponding property for ${\mathcal X}$-convex functions 
is the local Lipschitz continuity with respect to the Carnot-Carath\'eodory metric associated to the vector fields: this was proved in \cite{LMS, DGN, Mag, Ric, JLMS}  for the generators of a Carnot group  and in \cite{BD} for general 
 fields.

\smallskip
2. The operator in  \eqref{MAh} does not satisfy 
the standard structure conditions in viscosity theory, unless the vector fields are constant. To overcome this problem, for $H>0$ we take the $\log$ of both terms in \eqref{MAh} and show that the new equation verifies the Lipschitz-type condition 
with respect to $x$ of Crandall, Ishii, and Lions \cite{CIL} for {\em uniformly ${\mathcal X}$-convex} subsolutions, i.e., functions such that, for a $\gamma>0$, 
$
-D^2_{\mathcal X}u + \gamma I \leq 0
$
in $\Omega$ in viscosity sense. Our first main result states the comparison among semicontinuous sub- and supersolutions of 
equations of the form
\begin{equation}
\label{MAlog0}
-\log \det (D_{\mathcal X}^2 u)+ K(x, u, Du, D^2 u)=0, \ \mbox{in}\ \Omega,
\end{equation}
provided that either $K$ is strictly increasing in $u$ or that the subsolution is strict. Here $K$ is any degenerate elliptic operator satisfying the structure conditions of \cite{CIL}.

\smallskip
3. To cover the cases of $H$ not strictly increasing in $u$, which is the most frequent in applications, and satisfying only $H\geq 0$, we need to perturb a ${\mathcal X}$-convex subsolution to a uniformly ${\mathcal X}$-convex strict subsolution. This was done in the Euclidean case in \cite{IL} and we adapted the method to several nonlinear subelliptic equations in  \cite{BM}. We are able to perform this construction for equation \eqref{MAh} under an additional condition on the vector fields, namely
\begin{equation}
\label{C-type}
X_j(x) = \frac{\de}{\de x_j} + \sum_{i=m+1}^n \tau_{ij}(x) \frac{\de}{\de x_i} \quad j=1,...,m .
\end{equation}
In this case we say the  vector fields are of {\em Carnot type}, because this property is satisfied by the generators of a Carnot group. However we do not need all the other rich properties of such generators, not even the H\"ormander bracket generating condition.

 We therefore get the following Comparison Principle, under 
  the same assumptions as the Euclidean result of Ishii and 
  Lions \cite{IL}.

\begin{theorem}
\label{teo:main}
Assume $H : \ol\Omega\times\R\times \R^m \to [0, +\infty[$ is continuous, nondecreasing in the second entry, and for all $R>0$ there is $L_R$ such that
\begin{equation}
\label{H1/m}
| H^{1/m}(x,r,q+q_1) -  H^{1/m}(x,r,q) | \leq L_R|q_1| 
\quad\forall \; x\in\ol\Omega, |r|\leq R, |q|\leq R, |q_1|\leq 1.
\end{equation}
Suppose the vector fields $X_1,...,X_m\in C^2$ satisfy \eqref{C-type}. 
Let $u\,:\,\ol\Omega\to\R$ be a bounded, ${\mathcal X}$-convex, u.s.c. subsolution of \eqref{MAh} and $v\,:\,\ol\Omega\to\R$ be a bounded l.s.c. supersolution of \eqref{MAh}. Then
\begin{equation}
\label{comparison}
\sup_\Omega(u-v) \leq  \max_{\de\Omega}(u-v)^+.
\end{equation}
\end{theorem}
Note that the result applies to the prescribed horizontal Gauss curvature equation \eqref{MAh} \eqref{Gauss}.
We 
also get a Comparison Principle for 
\eqref{MAf} 
that extends to the viscosity context 
 a result of Rauch and Taylor \cite{RT} in $W^{2,n}$, the first for not necessarily convex supersolutions.

Theorem \ref{teo:main} implies that there is at most one  ${\mathcal X}$-convex continuous viscosity solution of the equation \eqref{MAh} with prescribed 
 boundary data
\begin{equation}
\label{BC}
u=g \quad\text{on }{\de\Omega},\quad g\in C(\de\Omega).
\end{equation}
The existence of solutions to the Dirichlet problem can be studied by the Perron method, as adapted to viscosity solutions by Ishii \cite{CIL}. It turns out to fit very well with our modified notions of sub- and supersolution.
A byproduct is the following subelliptic version of a classical result by Caffarelli, Nirenberg, and Spruck \cite{CNS} and P.-L. Lions \cite{PLL2} in the Euclidean setting. 
\begin{theorem}
\label{esistenza1}
Under the assumptions of Theorem \ref{teo:main} suppose that $g\in C^2(\ol\Omega)$ is ${\mathcal X}$-convex and $\O$ is ${\mathcal X}$-convex, i.e., it is the sublevel set of a $C^2$ ${\mathcal X}$-convex function. Then the solvability of \eqref{MAh} \eqref{BC} is equivalent to the existence of a ${\mathcal X}$-convex subsolution attaining continuously the boundary data.
\end{theorem}
The construction of a subsolution with the desired properties requires further assumptions, as it is well-known in the Euclidean case \cite{GT, PLL2}. Our main existence result is the following extension of a Theorem of Lions \cite{PLL2}.
\begin{theorem}
\label{esistenza2}
Besides the assumptions of Theorem \ref{teo:main} suppose 
 \begin{equation}
 \label{crescita0}
H^{1/m}\left(x,\max_{\de\O}g,p\right)\leq L|p| + M, \quad \forall\; x\in\ol\O,\; p\in\R^m .
\end{equation}
Assume also that 
$\O$ is uniformly ${\mathcal X}$-convex, i.e., the sublevel set of a uniformly $C^2$ ${\mathcal X}$-convex function. Then there is a unique ${\mathcal X}$-convex solution in $C(\ol\O)$ of the Dirichlet problem \eqref{MAh} \eqref{BC}.
\end{theorem}
The growth condition \eqref{crescita0} rules out the prescribed Gauss curvature equation \eqref{Gauss}, where it is known that $k$ must satisfy some compatibility conditions \cite{GT, PLL2}. We have an existence result in this case for the Koranyi ball of the Heisenberg group if $k(x)\leq k_{\mathbb H}(x)$, where $k_{\mathbb H}$ is the horizontal Gauss curvature of the graph of the gauge $w$, i.e., $w(x)=|x|_{\mathbb H}^4$ and $|x|_{\mathbb H}$ is the homogeneous norm of the Heisenberg group.

Some of the comparison results proved here were announced in 
our papers \cite{BM2} and \cite{BM4}. More precisely, \cite{BM2} contains the statement of Theorem \ref{teo:main} in the case of Carnot groups and $H>0$, with 
some ideas of the proof, and \cite{BM4} gives a different proof of Theorem \ref{confrontoC}  below for strictly increasing $H$ and in Carnot groups.

The Dirichlet problem for  subelliptic fully nonlinear equations was studied by Bieske \cite{B1, B2}, Bieske and Capogna \cite{BC}, and Wang \cite{W3} for the Aronsson equations of the Calculus of Variations in $L^\infty$, and by ourselves \cite{BM} and Cutri and Tchou \cite{CT} for Pucci-type and other Bellman-Isaacs equations. 
Almost nothing is known on the regularity of solutions of fully nonlinear subelliptic equations. This is a challenging subject for future research.

The paper is organized as follows. Section \ref{CPstrict} is devoted to 
the definitions, the Comparison Principle for the equation \eqref{MAlog0} (and variants of it), and its applications to \eqref{MAh} if $H>0$ and either the subsolution is strict or $H$ is strictly increasing in $r$. In Section \ref{CPgeneral}
we build strict subsolutions for vector fields of Carnot type and complete the proof of Theorem \ref{teo:main}. 
Section \ref{DP}  deals with the existence issue for the Dirichlet problem.

\section{Comparison principles with strict subsolutions}\label{CPstrict}
\subsection{Definitions} 

Let us consider 
equations of the form
\begin{equation}
\label{MAsube}
-G
 \left(\sigma^T(x) D^2 u\, \sigma(x)+A(x,Du)\right)+ K(x, u, Du, D^2u)=0 \;\; \mbox{in}\ \Omega,\quad G=\det \text{ or }G=\log \det ,
\end{equation}
where the set $\Omega\subseteq\re^n$ is open  and bounded.
\noindent We denote with $S^n$  the set of the symmetric $n\times n$ matrices, with $\leq$ the  usual partial order, 
with $I$ the identity matrix, and with $\mbox{tr} M$ the trace of a square matrix $M$. 
By $M>0$ we denote any positive definite matrix.
$USC(\ol\Omega)$ and $LSC(\ol\Omega)$ denote the sets of functions $\ol\Omega\to\R$ that are, respectively, upper semicontinuous and lower semicontinuous.
The assumptions on the data are the following.
\begin{eqnarray}\label{B}
&&K: \ol \Omega\times \re\times \re^n\times S^n   \rightarrow \re\ 
\mbox { is continuous};\\
&&K(x, r, p,X)\leq K(x, s, p, Y),\ \forall r\leq s,\ Y\leq X\nn\\
&&\forall x\in \ol\Omega, u\in\re, p\in\re^n, X,Y\in S^n.\nn
\end{eqnarray}
\begin{eqnarray}\label{STRUT}
&&K\left(y, r, \frac {x-y}{\epsilon},Y\right) - K\left(x, r,  \frac {x-y}{\epsilon},X\right)\leq \omega\left( |x-y|\left(1+\frac {|x-y|}{\epsilon}\right)\right),\\
&&\mbox{ for some modulus } \omega \mbox{ and all } \epsilon>0, x, y\in \ol\Omega,\ r\in \R, X,Y\in S^n \mbox{ satisfying } \nn\\
&& \ -\frac{3}{\epsilon}\left(\begin{array}{cc}
I&0\\
0&I\end{array}\right)\leq
\left(\begin{array}{cc}
X&0\\
0&-Y\end{array}\right)\leq
\frac{3}{\epsilon}\left(\begin{array}{cc}
I&-I\\
-I&I\end{array}\right),\label{structure}
\end{eqnarray}
\begin{eqnarray}\label{sigma}
&&\sigma(x) \mbox { is a Lipschitz continuous }  n\times m \mbox { matrix  valued function on } \ol\Omega, m\leq n. 
\end{eqnarray}
\begin{eqnarray}\label{Q0}
&& A(x,p) \mbox { is a continuous }  m\times m \mbox { matrix valued function on } 
\ol\Omega\times \re^n
\mbox { such that }\\
&&-C_1 |x-y|(1+|p|)I\leq A(x,p)-A(y,p)\leq C_1 |x-y|(1+|p|)I.\nn
\end{eqnarray}
\vskip 0.5cm

%

\begin{definition} 
\label{matrixineq}
If $\Psi: \ol \Omega\times \re^n\times S^n   \rightarrow S^m$ and $M\in S^m$
we say that $u$ is a (viscosity) subsolution of the matrix inequality  
$\Psi(x, Du, D^2u)\leq M ,\ in\ \Omega,
$
 if $u$ is USC in $\Omega$ and 
$
\Psi(x, D\phi(x), D^2\phi(x))\leq M ,
$  
 for all $\phi\in C^2(\Omega)$ and $x\in\argmax(u-\phi)$.
\end{definition}
The definition of (viscosity) subsolution $u$ of (\ref{MAsube}) is given in a standard way, as in \cite{CIL} (see also the comments in the next Remark \ref{commentidef}).
 \begin{definition}
 \label{defsubsol}
A function  $u\in USC(\ol\Omega)$
 is a (viscosity) subsolution  of  \eqref{MAsube} with $G= \det$ or $G=\log\det$ if 
for all  $\phi \in C^2(\ol\Omega)$ such that $u-\phi$ has a
 maximum point at $x_0$  we have
\begin{equation}\label{MAsubsolu}
- G
 (\sigma^T(x_0) D^2 \phi(x_0)\, \sigma(x_0)+A(x_0, D\phi(x_0))
+ K(x_0, u(x_0), D\phi(x_0), D^2\phi(x_0))\leq 0. 
\end{equation}
\end{definition}
The definition of (viscosity) supersolution $v$ of (\ref{MAsube}) is modified 
by restricting the test functions 
to the $C^2$ functions $\phi$ with 
$\sigma^T
D^2 \phi
\, \sigma
+A
> 0$
 at points $x\in\arg\min (v-\phi)$. The motivation is the same as in the Euclidean case \cite{IL} where, however,
 test functions with $D^2 \phi\geq 0$ were allowed. (See also \cite{Tru90} for viscosity solutions of other prescribed curvature equations).
\begin{definition} 
\label{defsupersol}
A function
$v\in LSC(\ol\Omega)$
 is a (viscosity) supersolution 
of  \eqref{MAsube} if 
for all  $\phi \in C^2(\ol\Omega)$ such that $v-\phi$ has a
 minimum point at $x_0$ and 
 \begin{equation}\label{stella}
 \sigma^T(x_0) D^2 \phi(x_0)\, \sigma(x_0)+A(x_0, D\phi(x_0))> 0, 
 \end{equation}
 we have, for $G= \det$ or $G=\log\det$,
\begin{equation}\label{MAsupersolu}
-G
 (\sigma^T(x_0) D^2 \phi(x_0)\, \sigma(x_0)+A(x_0, D\phi(x_0))
+ K(x_0, v(x_0), D\phi(x_0), D^2\phi(x_0))\geq 0.
\end{equation}
\end{definition}

\begin{remark}
\label{commentidef}
{\rm
In the next sections we will compare a supersolution of \eqref{MAsube} in the sense of Definition \ref{defsupersol} with a function $u$ subsolution of \eqref{MAsube} as in Definition \ref{defsubsol} satisfying also the matrix inequality 
\begin{equation}\label{matrix}
-\left(\sigma^T(x) D^2 u\, \sigma(x) + A(x,Du)\right)\leq 0 
\end{equation}
 in the sense of Definition \ref{matrixineq}. This is equivalent to comparing sub- and supersolutions in the standard sense of  \cite{CIL} of the equation
$$
\max\left\{-\lambda_{min}\left(\sigma^T D^2 u\, \sigma + A(x,Du)\right), -G
 \left(\sigma^T D^2 u\, \sigma+A(x,Du)\right)+ K(x, u, Du, D^2u)\right\} =0 ,
$$
where $\lambda_{min}(Z)$ denotes the minimal eigenvalue of $Z\in S^m$.
This is obvious for subsolutions, whereas  for a supersolution $v$ of the last equation and a standard test function $\phi$, either \eqref{MAsupersolu} holds, or $\lambda_{min}\left(\sigma^T(x_0) D^2 \phi(x_0)\, \sigma(x_0)+A(x_0, D\phi(x_0)\right)\leq 0$ at $x_0\in\arg\min (v-\phi)$, which is equivalent to Definition \ref{defsupersol}.

In the case $G= \log\det$ we will further restrict subsolutions to functions satisfying the matrix inequality 
$-\left(\sigma^T(x) D^2 u\, \sigma(x) + A(x,Du)\right)\leq -\gamma I$  in the sense of Definition \ref{matrixineq}. Then the first term in \eqref{MAsubsolu} is well defined because the argument of $G$ is a positive definite matrix.
}\end{remark}

The generality of the term $A$ in \eqref{MAsube} includes equation of the form 
$$
-\log\det \left(D^2u + A(x)\right)+f(x,u)=0 ,
$$
arising in problems of Riemannian geometry (see \cite{CNS, Au} and the references therein) and their counterparts involving non-commutative vector fields. However, in this paper we are mostly interested in 
subelliptic equations
\begin{equation}
\label{MAsubeN}
-\det D_{\mathcal X}^2 u+ F(x, u, D_ {\mathcal X} u, D_{\mathcal X}^2 u)=0, \quad \mbox{in}\ \Omega,
\end{equation}
where 
$D_ {\mathcal X}u = (X_1u,...,X_mu)$ is the intrinsic (or horizontal) gradient
with respect of a family of $C^{1,1}$ 
vector fields $X_1,...,X_m$,  and 
$$(D_ {\mathcal X}^2 u)_{ij} = \frac{X_i(X_ju)+ X_j(X_iu)}{2}$$
 is the symmetrized intrinsic Hessian.
 If we take the $n\times m$   $C^{1,1}$ matrix-valued function $\sigma$,
defined in $\ol\Omega\subseteq \R^n$,  whose columns $\sigma^j$ are the coefficients 
 of $X_j$, $j=1,\cdots, m$, we see that, for any smooth $u$
\begin{equation}
\label{completo0}
D_ {\mathcal X}u = \sigma^T (x)Du, \quad  D_ {\mathcal X}^2 u= \sigma^T(x) D^2 u\, \sigma(x) + Q(x, Du), 
\end{equation}
where 
$Q(x,p)$ is a $m\times m$ matrix whose elements are
\begin{equation}\label{Q}
Q_{ij}(x, p):= \frac12\left(
D\sigma^j(x)\,\sigma^i(x)+D\sigma^i(x)\,\sigma^j(x)
\right)\cdot p .
\end{equation}
Therefore the PDE \eqref{MAsubeN} can be written in the form \eqref{MAsube} with 
$$A=Q,\quad
G=\det ,\quad
K(x,r,p,X)=F(x,r,\s^Tp,\sigma^T X\, \sigma +Q).
$$
In this case the functions satisfying the matrix inequality \eqref{matrix} are called $\mX$-convex, consistently with the theory of convex functions in Carnot groups \cite{LMS, JLMS} and in general Carnot-Caratheodory metric spaces \cite{BD}.
\begin{definition}\label{Xconv}
$u\in USC(\ol\Omega)$ is convex  in $\Omega$ with respect to the fields $X_1,\dots,X_m$, briefly $\mathcal X$-convex (resp., uniformly $\mathcal X$-convex), if it is a subsolution of  
\begin{equation}\label{hconvex}
- D_ {\mathcal X}^2 u=-\sigma^T(x) D^2 u\, \sigma(x)-Q(x,Du)\leq 0,\ \mbox{in}\ \Omega
\end{equation}
(resp., $\leq-\gamma I$ for some $\gamma>0$).
\end{definition}

\begin{remark}{\rm
Note that for a uniformly $\mathcal X$-convex subsolution of \eqref{MAsubeN} the 
test functions can be restricted to $C^2$ strictly $\mX$-convex functions (i.e., satisfying \eqref{stella}), as for supersolutions.
}
\end{remark}

\subsection{The basic Comparison Principle}

The first result is a Comparison Principle between a supersolution and a strict subsolution such that 
$-\sigma^T D^2 u\, \sigma - A(x,Du)\leq -\gamma\, I$ of the equation
\begin{equation}
\label{MAlog}
-\log \det
 \left(\sigma^T(x) D^2 u\, \sigma(x)+A(x,Du)\right)+ K(x, u, Du, D^2u)=0 \;\; \mbox{in}\ \Omega .
\end{equation}

\begin{theorem}
\label{confronto}
Assume \eqref{B},  \eqref{STRUT}-\eqref{structure}, \eqref{sigma},  \eqref{Q0}.
Let $u\in USC(\ol\Omega)$ be a bounded subsolution, for some $\gamma, \gamma_1>0$,  of
\begin{equation}\label{conv}
-\sigma^T(x) D^2 u\, \sigma(x) -A(x,Du)\leq -\gamma\, I ,\ in\ \Omega,
\end{equation}
and 
\begin{equation}\label{strictPDE}
-\log\det (\sigma^T(x) D^2 u\, \sigma(x)+A(x, Du))+ K(x, u, Du, D^2u)\leq -\gamma_1,\ in\ \Omega.
\end{equation}
Let $v\in LSC(\ol\Omega)$ be a bounded supersolution of \eqref{MAlog}. Then
\begin{equation}\label{cp}
\sup_{\Omega}(u-v)\leq \max_{\partial \Omega}(u-v)^+.
\end{equation}
\end{theorem}
To prove the Comparison Principle we need the following two lemmata.
\begin{lemma}\label{rappres}
If $\gamma>0$, for all $A\in S^N$, $A\geq \gamma\, I$,
\begin{eqnarray}
&&\log\det(A)=\label{represlog}\\
&&=\min\{ N\log a-N+\mbox{tr}(AM) : a>0, M\in S^N, 0\leq M\leq \frac{1}{\gamma}I, \text{det}\ M=a^{-N}\}.\nn
\end{eqnarray} 
\end{lemma}
\begin{proof}
It is well known that
\begin{equation}\label{represdet}
(\det A)^{1/N}=\min\{ \mbox{tr}(AB), B\in S^N, B\geq 0, \det\, B=N^{-N}\},
\end{equation} 
and the minimum is attained at 
$B_m=\frac{(\det A)^{1/N}}{N}\,A^{-1}$.
On the other hand 
$$
\log[(\det A)^{1/N}]= \min\{ \log a+\frac{(\det A)^{1/N}-a}{a}: a>0\},
$$
and the minimum is attained at $a_m= (\det A)^{1/N}$.
We combine the two formulas to get a minimum representation for 
$\log[(\det A)^{1/N}]$ and we can restrict the search for the minimum  to matrices of the form
$B= \frac{a}{N} M$ with $M^{-1}\geq \gamma I$. Then 
\begin{eqnarray*}
&&\frac{1}{N}\log\det (A)=\\
&&= \min\{ \log a-1+\frac{\mbox{tr}(A\frac{a}{N}M)}{a} : a>0,\  M\in S^N, \ 0\leq M\leq \frac{1}{\gamma}I,\ a^{N}\det\, M=1\},
\end{eqnarray*}
which gives \eqref{represlog}.
\end{proof}

\begin{lemma}\label{lemmadiff}
Consider the operator
$$
F(x,p,X):= -\log \det ( \sigma^T(x) X\, \sigma(x)+A(x,p))$$
with $x\in \ol\Omega$, $p\in \R^n$, $X\in S^n$, $\sigma$  and $A$ satisfying  \eqref{sigma},  \eqref{Q0}.
Then for all $\gamma>0$ 
there is a constant $C>0$ such that
\begin{eqnarray}
&&F\left(y,\frac{x-y}{\epsilon}, Y\right)-F\left(x,\frac{x-y}{\epsilon}, X\right)\leq C\left(|x-y|+\frac{|x-y|^2}{\epsilon}\right)\label{diff}
\end{eqnarray}
 for all  $X, Y \in S^n$ 
 satisfying \eqref{structure} and 
\begin{equation}
 \sigma^T(x) X\, \sigma(x)+A\left(x,\frac{x-y}{\epsilon}\right)\geq \gamma I,\  \sigma^T(y) Y\, \sigma(y)+A\left(y,\frac{x-y}{\epsilon}\right)\geq \gamma I.\label{gamma}
\end{equation}
\end{lemma}
\begin{proof}
By Lemma \eqref{rappres} and \eqref{gamma} we can write $F(x,p,X)$ as the maximum of 
$$m-m\log a-\text{tr} (\sigma^T(x) X\, \sigma(x)M)-\text{tr} (A(x,p)M)$$
 as $a$, $M$ vary over $a>0$, $M\in S^m$, $0\leq M\leq \frac{1}{\gamma}I$,
$\det\, M=a^{-m}$.
Then there is a choice of $a$ and $M$ such that the left hand side of (\ref{diff}) is bounded above by the sum 
of
\begin{equation}\label{1tr}
\text{tr} (\sigma^T(x) X\, \sigma(x)M)-\text{tr} (\sigma^T(y) Y\, \sigma(y)M),
\end{equation}
and
\begin{equation}\label{2}
\text{tr} (A(x,p)M)- \text{tr} (A(y,p)M),\quad p=\frac{x-y}{\epsilon}.
\end{equation}
By diagonalization we see that $|M|\leq \frac{\sqrt m}{\gamma}$, where  $|\cdot|$ denotes the Euclidean norm.
Moreover there is $R\in S^m$ such that $M=R R^T$, $|R|\leq \sqrt{\frac{\sqrt m}{\gamma}}$.
We call $\Sigma(x)$ the $n\times n$ matrix whose first lines are $R\sigma^T(x)$ and the last $m-n$ lines are $0$.
Then (\ref{1tr}) can be rewritten as 
$$
\text{tr} (\Sigma^T(x) \Sigma (x) X)-\text{tr} (\Sigma^T(y) \Sigma (y) Y).
$$
A standard calculation in the theory of viscosity solutions (see, e.g., Example 3.6 in \cite{CIL}) shows that 
this quantity is bounded above by $ 3L^2\frac{|x-y|^2}{\epsilon}$ for matrices satisfying \eqref{structure}, where
$L$ is a Lipschitz constant of $\Sigma(\cdot)$. Therefore we can take $L=L_{\sigma} \sqrt{\frac{m}{\gamma}}$
where $L_{\sigma}$ is a Lipschitz constant for $\sigma(\cdot)$.
As for (\ref{2}),
$$
|\text{tr} \big(A(x,p)-A(y,p)\big)M|\leq |A(x,p)-A(y,p)||M|\leq  \frac{C_1\sqrt m}{\gamma}|x-y|(1+|p|).
$$
In conclusion, we get 
\begin{eqnarray*}
F(y,\frac{x-y}{\epsilon}, Y)-F(x,\frac{x-y}{\epsilon}, X)\leq \frac{C_1\sqrt m}{\gamma}|x-y|+
\frac{3L_{\sigma}^2m+C_1\sqrt m}{\gamma}\frac{|x-y|^2}{\epsilon}
\end{eqnarray*}
\end{proof}
We can prove now the comparison theorem.
\begin{proof} (of Theorem \ref{confronto})
For $\epsilon>0$  the function 
$\Phi_{\epsilon}(x,y)=u(x)-v(y)-\frac{1}{2\epsilon}|x-y|^2$
has a maximum point 
$(x_{\epsilon}, y_{\epsilon})$.
A standard argument gives 
\begin{equation}\label{epsto0}
\frac{|x_{\epsilon}- y_{\epsilon}|^2}{\epsilon}\rightarrow 0, \mbox{ as }  \epsilon \rightarrow 0^+.
\end{equation}
If there is a sequence $\epsilon_k \rightarrow 0$ such that $x_{\epsilon_k}\rightarrow \hat x\in\partial\Omega$,
then $y_{\epsilon_k}\rightarrow \hat x$, and by the upper semicontinuity of $u(x)-v(y)$
$$\max_{\ol\Omega}(u-v)\leq \Phi_{\epsilon}(x_{\epsilon}, y_{\epsilon})\rightarrow \max_{\partial\Omega}(u-v),\ \mbox{ as }
\epsilon \rightarrow 0.$$ 
The case of 
$y_{\epsilon_j}\rightarrow \hat y\in \partial\Omega$ for some $\epsilon_j \rightarrow 0$ is analogous. 
Therefore we are left with the case $(x_{\epsilon}, y_{\epsilon})\in\Omega\times\Omega$ for all small $\epsilon$.
We use the Theorem on Sums, as in \cite{CCime} and get $X,Y\in S^n$ (depending on $\e$) such that, for 
$p_{\epsilon}:=\frac{|x_{\epsilon}- y_{\epsilon}|}{\epsilon}$,
$(u(x_{\epsilon}), p_{\epsilon}, X)\in \ol J^{2,+}u(x_{\epsilon})$, 
$(v(y_{\epsilon}), p_{\epsilon}, Y)\in \ol J^{2,-}v(y_{\epsilon})$, 

\begin{eqnarray}
&& \ -\frac{3}{\epsilon}\left(\begin{array}{cc}
I&0\\
0&I\end{array}\right)\leq
\left(\begin{array}{cc}
X&0\\
0&-Y\end{array}\right)\leq
\frac{3}{\epsilon}\left(\begin{array}{cc}
I&-I\\
-I&I\end{array}\right).\label{MI}
\end{eqnarray}
Then (\ref{conv}) implies 
$$G(x_{\epsilon}, X):= \sigma^T(x_{\epsilon}) X\, \sigma(x_{\epsilon})+ A(x_{\epsilon}, p_{\epsilon})\geq \gamma I.$$
We seek a similar inequality for 
$G(y_{\epsilon}, Y):= \sigma^T(y_{\epsilon}) Y\, \sigma(y_{\epsilon})+ A(y_{\epsilon}, p_{\epsilon})$.
To this end we multiply on the left the second inequality in (\ref{MI}) by the 
 $m\times 2n$ matrix whose first $n$ columns are $\sigma^T(x_{\epsilon})$ and the last $n$ are 
 $\sigma^T(y_{\epsilon})$, and then on the right by the transpose of such matrix. Since the operation preserves the inequality, we get
\begin{eqnarray}
&&\sigma^T(x_{\epsilon}) X \sigma(x_{\epsilon})-\sigma^T(y_{\epsilon}) Y \sigma(y_{\epsilon}) \leq \frac{3}{\epsilon} 
(\sigma(x_{\epsilon})-  \sigma(y_{\epsilon}))^T(\sigma(x_{\epsilon})-  \sigma(y_{\epsilon}))\label{QI}\\
&&
\leq 
\frac{3}{\epsilon}C_{\sigma}|x_{\epsilon}-y_{\epsilon}|^2\,I,\nn
 \end{eqnarray}
 where $C_{\sigma}$ is a suitable constant related to the Lipschitz constant of $\sigma$. 
 Then, by (\ref{QI}) and assumptions (\ref{Q0}),
 \begin{eqnarray*}
&&G(y_{\epsilon}, Y)\geq   G(x_{\epsilon}, X)
- \frac{3}{\epsilon}C_{\sigma}|x_{\epsilon}-y_{\epsilon}|^2\,I + A(y_{\epsilon}, p_{\epsilon})- A(x_{\epsilon}, p_{\epsilon})\\
&&\geq \big(\gamma- \frac{3}{\epsilon}C_{\sigma}|x_{\epsilon}-y_{\epsilon}|^2-C_1|x_{\epsilon}-y_{\epsilon}|- 
C_1\frac{|x_{\epsilon}-y_{\epsilon}|^2}{\epsilon}\big)\,I \geq \frac{\gamma}{2}\, I\nn
 \end{eqnarray*}
 for $\epsilon$ small enough, by (\ref{epsto0}). 
 Now we use the fact that $u$ satisfies (\ref{strictPDE}) and $v$ is a supersolution 
 to get
 \begin{eqnarray}
&& -\log\text{det}(\sigma^T(x_{\epsilon}) X\, \sigma(x_{\epsilon})+ A(x_{\epsilon}, p_{\epsilon}))+ 
K(x_{\epsilon}, u(x_{\epsilon}), p_{\epsilon}, X)\leq -\gamma_1<0,\label{1}\\
 &&
 -\log\text{det}(\sigma^T(y_{\epsilon}) Y\, \sigma(y_{\epsilon})+ A(y_{\epsilon}, p_{\epsilon}))+ 
K(y_{\epsilon}, v(y_{\epsilon}), p_{\epsilon}, Y)\geq 0.\nn
  \end{eqnarray}
 If for some $\epsilon$, $u(x_{\epsilon})\leq v(y_{\epsilon})$ we conclude that 
 $\max_{\ol\Omega}(u-v)\leq u(x_{\epsilon})-v(y_{\epsilon})\leq 0$. Otherwise, by the monotonicity 
 of $K$ with respect to the 2nd entry $r$, we get
$$
  -\log\text{det}(\sigma^T(y_{\epsilon}) Y\, \sigma(y_{\epsilon})+ A(y_{\epsilon}, p_{\epsilon}))+ 
K(y_{\epsilon}, u(x_{\epsilon}), p_{\epsilon}, Y)\geq 0.$$
 Now we subtract this inequality from the first of (\ref{1}), we use Lemma \ref{lemmadiff} 
 and the structure condition on $K$ to obtain 
$$C|x_{\epsilon}-y_{\epsilon}|(1+ \frac{|x_{\epsilon}-y_{\epsilon}|}{\epsilon})+
\omega(|x_{\epsilon}-y_{\epsilon}|(1+ \frac{|x_{\epsilon}-y_{\epsilon}|}{\epsilon}))\leq-\gamma_1<0$$
which gives a contradiction as 
 $\epsilon\rightarrow 0^+$ , by (\ref{epsto0}).
\end{proof}

\begin{remark}
\label{no+}
{\rm
If $K=K(x,p,X)$ is independent of $r$ the previous proof shows that 
$$\sup_{\Omega}(u-v)\leq \max_{\partial \Omega}(u-v).$$
}
\end{remark}

\begin{remark}\label{confronto1}
{\rm
Theorem \ref{confronto} remains true if we relax the strict subsolution condition (\ref{strictPDE}) to the following: 
for any open $\Omega_1$ such that $\ol\Omega_1\subseteq \Omega$ there exists $\gamma_1>0$ such that
$$
-\log\det (\sigma^T(x) D^2 u\, \sigma(x)+A(x, Du))+ K(x, u, Du, D^2u)\leq -\gamma_1,\ \mbox{in}\ \Omega_1
$$
holds.
The proof is the same, because if no sequence $x_{\epsilon_k}$ or $y_{\epsilon_j}$ converge to a boundary point then 
$(x_{\epsilon},  y_{\epsilon})\in \ol\Omega_1\times \ol\Omega_1$ for some 
$\ol\Omega_1\subseteq \Omega$ and for all $\epsilon$ small enough.
}
\end{remark}

\begin{theorem} 
\label{confrontostretto}
The conclusion of Theorem 
\ref{confronto} remains true if $u$ is a  subsolution of \eqref{MAlog} 
instead of a strict subsolution \eqref{strictPDE},
provided that, for some $C>0$,
\begin{eqnarray*}
&&K(x,r,p,X)-K(x,s,p,X)\geq C(r-s),\quad -M\leq s\leq r\leq M,\\
&& M:=\max 
\{\|u\|_\infty, \|v\|_\infty\},\ x\in\ol\Omega,\ p\in\re^n,\ X\in S^n.
\end{eqnarray*}
Under this condition there is at most one  viscosity solution $u$ of \eqref{MAlog}
such that $-\sigma^T(x) D^2 u\, \sigma(x)-A(x, Du)\leq -\gamma\, I$
 with prescribed continuous boundary data.
\end{theorem}
\begin{proof}
It is a standard variant of the preceding one.
\end{proof}




\begin{remark}{\rm
Note that, if we consider equation without $\log$, the structure condition (\ref{diff}) 
can be not  true. 
Take for example $\hat F(x,p,X):= -\det(\sigma^T(x) X\, \sigma(x)+A(x,p))$ with 
$A\equiv 0$ and 
$$
\sigma=\left(\begin{array}{cc}1 & 0 \\0 & 1 \\2y & -2x\end{array}\right),
$$
 the matrix associated to the Heisenberg group.
If $X$ and $Y$ satisfy the matrix inequality (\ref{structure}) then (\ref{QI}) holds.
If we take $Y$ diagonal it is easy to see that
$$\det \left(G(y_{\epsilon}, Y)+ \lambda I \right)= \det G(y_{\epsilon}, Y) +\lambda ^2 +\lambda \mbox{tr } G(y_{\epsilon}, Y),$$
where
$G(y_{\epsilon}, Y):= \sigma^T(y_{\epsilon}) Y\, \sigma(y_{\epsilon})$. 
Then from (\ref{QI}), taking $\lambda =3C_{\sigma}\frac{|x_{\epsilon}-y_{\epsilon}|^2}{\epsilon}$,
$$\det G(y_{\epsilon}, Y)-\det G(x_{\epsilon}, X)\geq 
\det G(y_{\epsilon}, Y)-
\det \left(G(y_{\epsilon}, Y)+ \lambda I \right)\geq -\lambda ^2-\lambda \mbox{tr }G(y_{\epsilon}, Y).$$
The term $\lambda \mbox{tr }G(y_{\epsilon}, Y)$ does not necessarily tend to zero as $\epsilon\rightarrow 0$. 
Taking, for example, $v(y)$ Lipschitz continuous, since 
the function 
$\Phi_{\epsilon}(x,y)=u(x)-v(y)-\frac{1}{2\epsilon}|x-y|^2$, introduced in the proof of Theorem \ref{confronto},
has a maximum point in 
$(x_{\epsilon}, y_{\epsilon})$, 
$\Phi_{\epsilon}(x_{\epsilon}, y_{\epsilon})\geq  \Phi_{\epsilon}(x_{\epsilon}, x_{\epsilon})$, i.e.
$|v(x_{\epsilon}) -v(y_{\epsilon})|\geq \frac{1}{2\epsilon}|x_{\epsilon} - y_{\epsilon}|^2$ and from the Lipschitz continuity 
of $v$ with constant $L$, we obtain
$|x_{\epsilon} - y_{\epsilon}|\leq 2\epsilon L$. In this case, since 
$\mbox{tr }G(y_{\epsilon}, Y)\leq \frac{C}{\epsilon}$ for $C>0$,  the best one can say is that,  for a suitable 
$K>0$,
$\lambda \mbox{tr }G(y_{\epsilon}, Y)\leq  K\frac{|x_{\epsilon}-y_{\epsilon}|^2}{\epsilon ^2}\leq 4 K L^2$.
} 
\end{remark}


\subsection{Subelliptic equations under a boundedness condition}\label{casosube}

In this section we apply the basic Comparison Principle to equations of the form
\begin{equation}
\label{MAsube2}
-\det D_{\mathcal X}^2 u+ F(x, u, D_ {\mathcal X} u, D_{\mathcal X}^2 u)=0, \quad \mbox{in}\ \Omega,
\end{equation}
where 
$D_ {\mathcal X}u = (X_1u,...,X_mu)$ is the intrinsic gradient
with respect of the vector fields $X_1,...,X_m\in C^{1,1}$, 
$(D_ {\mathcal X}^2 u)_{ij} = 
\left(X_i(X_ju)+ X_j(X_iu)\right)/{2}$ is the symmetrized intrinsic Hessian.
 If $\sigma$
 is the $n\times m$   $C^{1,1}$ matrix-valued function defined in $\ol\Omega\subseteq \R^n$ 
  whose columns $\sigma^j$ are the coefficients 
 of $X_j$, $j=1,\cdots, m$, 
 then for smooth $u$
\begin{equation}\label{completo}
D_ {\mathcal X}u = \sigma^T (x)Du, \quad  D_ {\mathcal X}^2 u= \sigma^T(x) D^2 u\, \sigma(x) + Q(x, Du), 
\end{equation}
where the matrix
$Q(x,p)$ is given by \eqref{Q}. 
%

\begin{lemma}\label{strictlogG}
Let $0<F(x,u,p,X)\leq C_1$, for any $ x\in \ol\Omega, u\in\re, p\in\re^m, X\in S^m$.
Let $w\in USC(\ol\Omega)$ be an uniformly $\mathcal X$-convex bounded function satisfying
$$
- \det 
D_ {\mathcal X}^2 w+  F(x, w, D_ {\mathcal X} w, D_{\mathcal X}^2 w)\leq -\alpha <0,\ \mbox{in } \Omega.
$$
Then there is $\alpha_1>0$ such that
\begin{equation}\label{strictlog2G}
- \log \det (
D_ {\mathcal X}^2 w)+  \log F(x, w, D_ {\mathcal X} w, D_{\mathcal X}^2 w)
\leq -\alpha_1 <0,\ \mbox{in } \Omega.
\end{equation}
\end{lemma}
 \begin{proof} By the properties of log
\begin{multline*}
\log \big( F(x, w, D_ {\mathcal X} w, D_{\mathcal X}^2 w) +\alpha\big) =
\log F(x, w, D_ {\mathcal X} w, D_{\mathcal X}^2 w)
+\log(1+\frac{\alpha}{F(x, w, D_ {\mathcal X} w, D_{\mathcal X}^2 w)})
\\
\geq \log F(x, w, D_ {\mathcal X} w, D_{\mathcal X}^2 w) +\log(1+\frac{\alpha}{C_1}),
\end{multline*}
which gives (\ref{strictlog2G}) with $\alpha_1=\log(1+\frac{\alpha}{C_1})$.
\end{proof}

%

\begin{corollary}\label{cor1sub}
Suppose that $0<F(x,u,p,X)\leq C_1$, $K=\log F$ satisfies \eqref{B} 
\eqref{STRUT}, and 
the $X_1,...,X_m$ are $C^{1,1}$ 
in $\ol\Omega$. 
Then the Comparison Principle holds between a uniformly $\mathcal X$-convex strict subsolution $u$ and 
a supersolution $v$ of equation \eqref{MAsube2}. 
\end{corollary}
\begin{proof}
We apply Theorem \ref{confronto} with $K=\log F$ and $A=Q$ given by (\ref{Q}), noting that the strict subsolution of (\ref{MAsube2})
is a strict subsolution of  (\ref{MAlog}) 
by  Lemma \ref{strictlogG}.
\end{proof}%
From Theorem \ref{confrontostretto} we immediately get:
\begin{corollary}
\label{cor2sub}
The Comparison Principle 
 is true if 
$\sigma$ is a $C^{1,1}$ $n\times m$ matrix valued function,
$F>0$,  $K=\log F$ satisfies \eqref{B} 
and \eqref{STRUT},
 and $u$ is a  uniformly $\mathcal X$-convex subsolution of \eqref{MAsube2}, 
not necessarily strict, provided that, for some $C>0$,
$$
F(x,r,p,X)-F(x,s,p,X)\geq C(r-s), \quad -M\leq s\leq r\leq M,\; M:=\max 
\{\|u\|_\infty, \|v\|_\infty\}, 
$$
for all $x\in\ol\Omega,\ p\in\re^m, X\in S^m$.
\end{corollary}
Note that here the upper bound 
 for $F$ is not required.


 
\subsection{Monge-Amp\`ere equation with unbounded gradient terms}

In this section we prove a comparison result 
for the equation
\begin{equation}\label{equazioneCarnot}
- \log \det (D_ {\mathcal X}^2 u)+ K_1(x, u, D_ {\mathcal X}u)=0.
\end{equation}
with Hamiltonian $K_1$ unbounded but independent of the second derivatives.
In this case we can exploit the boundedness of $D_ {\mathcal X}u$ for any $\mathcal X$-convex function $u$ to decrease the assumptions on the Hamiltonian
$K(x, r, p)=K_1(x, r, \sigma^T(x) p)$. The gradient estimate in viscosity sense of the next proposition was proved very recently by the first author and Dragoni for general vector fields \cite{BD}.
In the special case of Carnot groups, various authors showed, under different assumptions, the Lipschitz
continuity of $\mathcal X$-convex
functions with respect to the intrinsic metric of the group and bounds on their horizontal gradient in the sense of distributions \cite{LMS}, 
\cite{DGN}, \cite{Mag}, \cite{Ric}, \cite{JLMS}. From those results one can obtain a short proof of the 
gradient bound in viscosity sense, that we give for the convenience of readers mostly interested in the Carnot group setting (see the next Section \ref{carnot} for the definitions).
 
 \begin{proposition}
 \label{stimagradiente}
  \cite{BD}
Let the vector fields $X_1,...,X_m$ be of class $C^2$ and $u$ be $\mX$-convex and bounded in $\ol\Omega$. 
Then, for every open $\Omega_1$
with $\ol\Omega_1\subseteq \Omega$, there exists a constant $C$ such that
$$ |\sigma^T(x)\, Du|\leq C,\  \mbox{in } \Omega_1
$$
in viscosity sense.
\end{proposition}
\begin{proof} {\it (in the case of Carnot groups)}
It is known that $\mathcal X$-convexity implies local Lipschitz continuity with respect to the Carnot-Caratheodory distance by a result of Magnani \cite{Mag} and Rickly \cite{Ric}, see also \cite{JLMS}.
In particular, $u$ is continuous in $\Omega$. We mollify $u$ by convolution with kernels adapted to the group structure, as in \cite{DGN, BLU}. The approximating $u_{\epsilon}$ converge to $u$ uniformly on compact subsets of $\Omega$, and they are smooth and $\mathcal X$-convex.
Moreover, from the proof of Theorem 9.1 of \cite{DGN} we get, for $R$ small enough,
$$
\sup_{B_C(x_0,R)}\,\left(\sum_{j=1}^m (X_j u_{\epsilon})^2\right)^{1/2} \leq \frac{2}{R}\sup_{B_C(x_0,3R)}\, |u|,
$$
where the balls $B_C$ are taken with respect to the gauge pseudo-distance and $X_j u$ denotes the derivative of $u$ along the  trajectory of the vector field $X_j$. Since $u_{\epsilon}$ is $C^{\infty}$, 
$X_j u_{\epsilon}(x)=\sigma^j(x)Du_{\epsilon}(x)$. Therefore there is a constant $C$ depending only on 
$\sup_{\ol\Omega}\, |u|$ and the pseudo-distance of $\Omega_1$ from $\partial \Omega$ such that
$$ |\sigma^T(x)\, Du_{\epsilon}|\leq C,\  \mbox{in } \Omega_1.
$$
By letting $\epsilon\rightarrow 0$, we obtain that $u$ is a viscosity subsolution of the same inequality.
\end{proof}

\begin{theorem}\label{confrontoC}
Assume 
$K_1:\ol\Omega \times\R\times\R^m \rightarrow \R$ is continuous and nondecreasing w.r.t. 
$r$
and $X_1,\cdots,X_m$ are of class $C^2$. 
Suppose $u\in USC (\ol\Omega)$ is bounded and uniformly ${\mathcal X}$-convex, i.e., $-D_ {\mathcal X}^2u +\gamma I\leq 0$, and a subsolution of 
(\ref{equazioneCarnot}), whereas $v\in LSC(\ol\Omega)$ is a bounded supersolution of (\ref{equazioneCarnot}). 
Finally, assume that either $K_1$
is strictly increasing in $r$, 
$K_1(x, r, q)- K_1(x, s, q)\geq C(r-s)$, for some $C>0$ and all $r,s\in[-M,M]$, $M=\max\{\|u\|_{\infty}, \|v\|_{\infty}\}$, or $u$ is a strict
subsolution of (\ref{equazioneCarnot}) in each $\ol\Omega_1$ with $\ol\Omega_1\subseteq \Omega$. Then
$$
\sup_{\Omega}(u-v)\leq \max_{\partial \Omega}(u-v)^+.
$$
\end{theorem}
\begin{proof}
We only show how we can avoid the structure condition (\ref{STRUT}) on the Hamiltonian in the proof of Theorem \ref{confronto}.
Since $p_{\epsilon}=\frac{x_{\epsilon}-y_{\epsilon}}{\epsilon}$ is in the superdifferential of $u$ at 
$x_{\epsilon}\in \Omega_1$, Proposition \ref{stimagradiente}  gives
$$|\sigma^T(x_{\epsilon})p_{\epsilon}|\leq C.$$
Moreover
$$|\sigma^T(x_{\epsilon})p_{\epsilon} -\sigma^T(y_{\epsilon})p_{\epsilon} |\leq 
L_{\sigma} \frac{|x_{\epsilon}-y_{\epsilon}|^2}{\epsilon}\rightarrow 0\ \mbox{ as } \epsilon \rightarrow 0,$$
where $L_{\sigma}$ is a Lipschitz constant of $\sigma$, and therefore, for $\epsilon$ small, 
$$|\sigma^T(y_{\epsilon})p_{\epsilon} |\leq C+1.$$
Let $\omega_1$ be the modulus of continuity of $K_1$ on 
$\ol\Omega\times [-M,M]\times \ol B(0,C+1)$. Then
\begin{eqnarray*}
&&|K( x_{\epsilon}, u(x_{\epsilon}), p_{\epsilon})-K( y_{\epsilon}, u(x_{\epsilon}), p_{\epsilon})|=\\
&&|K_1( x_{\epsilon}, u(x_{\epsilon}), \sigma^T(x_{\epsilon}) p_{\epsilon})-K_1( y_{\epsilon}, u(x_{\epsilon}), \sigma^T(y_{\epsilon})p_{\epsilon})|\\
&&\leq \omega_1(|x_{\epsilon}-y_{\epsilon}| +L_{\sigma}\frac{|x_{\epsilon}-y_{\epsilon}|^2}{\epsilon})
\rightarrow 0\ \mbox{ as } \epsilon \rightarrow 0.
\end{eqnarray*}
The rest of the proof is the same as that of Theorem \ref{confronto}, taking into account Remark \ref{confronto1}, and Theorem
\ref{confrontostretto}.
\end{proof}
A direct proof of this theorem under the strict monotonicity assumption on $K_1$ is given 
in our paper \cite{BM4}.
%



\section{The Comparison Principle for vector fields of Carnot type}
\label{CPgeneral}
\subsection{Carnot groups}
\label{carnot}
We begin with recalling some well-known definitions. We adopt the terminology and notations of the recent book \cite{BLU}.
Consider a group  operation $\circ$ on $\R^n=\R^{n_1}\times...\times\R^{n_r}$ with identity $0$, such that 
$$
(x,y)\mapsto y^{-1}\circ x 
\quad \text{ is smooth,}
$$ 
and the dilation $\d_\l\,:\,\R^n\to\R^n$ 
\[
\d_\l(x) = \d_\l(x^{(1)},...,x^{(r)}) := (\l x^{(1)}, \l^2 x^{(2)},..., \l^r x^{(r)}), \quad x^{(i)}\in \R^{n_i}.
\]
If $\d_\l$ is an automorphism of the group $(\R^n, \circ)$ for all $\l>0$, 
$(\R^n, \circ, \d_\l)$ is a homogeneous Lie group on $\R^n$. We say that $m=n_1$ smooth vector fields $X_1, ..., X_m$ on $\R^n$ generate $(\R^n, \circ, \d_\l)$, and that this is a (homogeneous) Carnot group, if
\begin{itemize}
\item $X_1, ..., X_m$ are invariant with respect to  the left translations on $\R^n$ $\tau_\a(x):= \a\circ x$ for all $\a\in\R^n$,
\item $X_i(0) = \de /\de x_i$, $i=1,...,m$,
\item the  rank  of the Lie algebra generated by $X_1, ..., X_m$ is $n$ at every point $x\in\R^n$.
\end{itemize}
We refer, e.g., to \cite{Bell, BLU} for the connections of this  definition with the classical one in the context of  abstract Lie groups and for the properties of the generators. We will use only  the following property, and refer to Remark 1.4.6, p. 59 of \cite{BLU} for more precise informations.

\begin{proposition} If  $X_1, ..., X_m$ are generators of a Carnot group, then
\[
X_j(x) = \frac{\de}{\de x_j} + \sum_{i=m+1}^n \s_{ij}(x) \frac{\de}{\de x_i}
\]
with $\s_{ij}(x)=\s_{ij}(x_1,...,x_{i-1})$ homogeneous polynomials of a degree $\leq n-m$.
\end{proposition}
The previous Proposition implies that
$\sigma(x)= \left(\begin{array}{cc}
I\\
\tau(x)\end{array}\right)$ where $I$ is the $m\times m$ identity matrix and $\tau(x)$ is a   $(n-m)\times m$ matrix. 

If $X_1,\cdots,X_m$ are the generators of a Carnot group $\mathcal{G}$, the definition (\ref{Xconv}) of 
$\mathcal X$-convexity 
coincides with the definition of
convexity in $\mathcal{G}$ in viscosity sense (v-convexity) of Lu, Manfredi, Stroffolini \cite{LMS}.
A more geometric notion of convexity in $\mathcal{G}$, called {\em horizontal convexity} (or weak H-convexity), was introduced and studied
in the same seminal paper \cite{LMS} and, independently, by Danielli, Garofalo, and Nhieu  \cite{DGN}.
The equivalence of the two notions was studied by several authors, first in the Heisenberg groups \cite{LMS},
\cite{BaRic}, and then in general Carnot groups \cite{W2}, \cite{Mag}, \cite{JLMS}.

\subsection{Construction of strict subsolutions}
 \label{PCC}

In this section we construct a uniformly ${\mathcal X}$-convex strict subsolution of the subelliptic Monge-Amp\`ere equation 
\begin{equation}
\label{MAsube0}
-\det D_{\mathcal X}^2 u+ H(x, u, D_ {\mathcal X} u)=0, \quad \mbox{in}\ \Omega,
\end{equation}
from a ${\mathcal X}$-convex  subsolution $u$ (here the $\mX$ derivatives are defined by \eqref{completo} and \eqref{Q}). 
We therefore get a Comparison Principle for usual viscosity subsolutions, without the strictness assumption.

The first assumption is on the vector fields and it is motivated by the properties of generators of Carnot groups recalled in the preceding section:
%
\begin{eqnarray}
&&\sigma(x) \mbox { is a }  n\times m \mbox { matrix such that }\label{sigma2C}\\
&&\sigma(x)= \left(\begin{array}{cc}
I\\
\tau(x)\end{array}\right) \mbox{ where } I \mbox{  is the } m\times m \mbox{ identity matrix and } \tau(x)\nn\\
 &&\mbox{ is a  } C^{1,1} (n-m)\times m
\mbox{ matrix. }\nn
\end{eqnarray}
  When the matrix $\sigma$ satisfies it we will say that the vector fields are 
  of {\em Carnot type}, following the terminology of \cite{MSC}. 
  However, 
  different from \cite{MSC}, we do not assume the H\"ormander condition on the rank of the Lie algebra generated by the fields.
  
 The second assumption is on $H$. In the Euclidean case it coincides with the one made by Ishii and P.-L. Lions for their Comparison Principle in that context \cite{IL}.
\begin{eqnarray}\label{LipCH}
&&H: \ol\Omega\times \R\times \R^m \rightarrow [0,+\infty),\\
&&H \mbox{ continuous and nondecreasing in } r,\nn\\
&&\mbox{for any } R >0\ \mbox {there exists}\ L_R\ \mbox {such that:}\nn\\
&&|H^{1/m}(x, r, q+q_1)- H^{1/m}(x, r, q)|\leq L_R|q_1|,\nn\\
&&\forall x\in \ol\Omega, |r|\leq R, |q|\leq R, |q_1|\leq 1.\nn
\end{eqnarray}
\begin{theorem}
\label{strict0}
Assume \eqref{sigma2C} and \eqref{LipCH} 
and let $u$ be a ${\mathcal X}$-convex 
subsolution 
of equation  \eqref{MAsube0}.
Then for any open set $\Omega_1$ with $\ol\Omega_1\subseteq \Omega$
there exist $\a, \e_o>0$ and a sequence $u_{\epsilon}\in USC(\ol\O)$ of uniformly ${\mathcal X}$-convex functions such that $u_{\epsilon}\leq u$, $u_{\epsilon} \rightarrow u$ uniformly in $\O$ as $\epsilon\to 0$, and 
\begin{equation}
\label{strict1}
-\det D^2_\mX u_\e 
+ H(x, u_{\epsilon}, D_\mX u_{\epsilon})\leq -\alpha, \quad \text{in } \Omega_1 ,\; \forall\,\, \e\leq\e_0 .
\end{equation}
\end{theorem}

\begin{proof}
We consider 
$$
u_{\epsilon}(x):=u(x)+\epsilon (e^{\mu\frac{\sum_{i=1}^m |x_i|^2}{2}}-\lambda),
$$
and we want to show that it is a strict subsolution 
for $\lambda$ and $\mu$ sufficiently large, 
independent of $\epsilon >0$. 
First we choose $\lambda:=\max_{x\in\ol\O}e^{\mu\frac{\sum_{i=1}^m |x_i|^2}{2}}$  for all $x\in\Omega$,
and this implies  $u_{\epsilon}(x)\leq u(x)$.
We set 
$$
\nu:= \epsilon \mu e^{\mu\frac{\sum_{i=1}^m |x_i|^2}{2}}, \qquad \epsilon_o:=\min_{x\in\ol\Omega}\frac{\exp (-\frac{\mu}{2}\sum_{i=1}^m\, x_i^2)}{\mu(\sum_{i=1}^m\, x_i^2)^{1/2}} ,
$$ 
and  compute
$$
D u_{\epsilon}=D u +\nu(x_1,\cdots,x_m, 0, \cdots,   0),$$
$$D ^2u_{\epsilon}=
D^2 u +\nu\left(\left(\begin{array}{cc}
I_m& 0\\
0 & 0\end{array}\right)+ \mu (x_1,\cdots,x_m, 0, \cdots,   0)\otimes(x_1,\cdots,x_m, 0, \cdots,   0)\right),
$$
where $(q\otimes q)_{ij}=q_i q_j$. Note that $|\nu(x_1,\cdots,x_m)|\leq 1$ for $\e\leq\e_o$.
Then
\begin{eqnarray*}
&&\sigma^T\, Du_{\epsilon}= \sigma^T Du +\nu(x_1,\cdots, x_m),\\
&&\sigma^T(x)\, D^2u_{\epsilon}\, \sigma(x) + Q(x, D u_{\epsilon} )
=\sigma^T(x)\, D^2u\, \sigma(x) +Q(x, D u)+\\
&&+\nu\left(I_m+\mu (x_1,\cdots, x_m)\otimes(x_1,\cdots, x_m)\right) +
\nu Q(x, (x_1,\dots,x_m, 0, \dots,   0)).
\end{eqnarray*}
From the structure of the coefficients of the matrices $Q$ and $\sigma$, \eqref{Q} and \eqref{sigma2C},
we have that 
$$
D\sigma^i= \left(\begin{array}{cc}
0\\
D\tau^i\end{array}\right),\qquad D\sigma^i\,\sigma^j= \left(\begin{array}{cc}
0\\
D\tau^i\,\tau^j\end{array}\right),
$$
 $\tau^i$ being the i-th column of the matrix $\tau$.
 Then
$Q_{ij}(x, (x_1,\cdots,x_m, 0, \cdots,   0))\equiv 0$ for any $i,j=1,\cdots m$.
Hence, since $u$ is ${\mathcal X}$-convex,
$$
-\sigma^T(x)\, D^2u_{\epsilon}\, \sigma(x) - Q(x, D u_{\epsilon})+\nu I_{m}\leq 0,
$$
 i.e., 
$u_{\epsilon}$ is uniformly ${\mathcal X}$-convex.

Now we want to find a sufficiently large $\mu$ such that 
$u_{\epsilon}$ satisfies (\ref{strict1}).
Let us consider the auxiliary equation
\begin{equation}\label{1/m}
G(x, u, Du, D^2u):= -\text{det}^{1/m} (\sigma^T\, D^2u\, \sigma+Q(x, D u) )+
H^{1/m}(x, u, \sigma^T Du)=0 .
\end{equation}
To prove that $u_{\epsilon}$  is a strict subsolution of (\ref{1/m}) for 
large $\mu$ we compute 
\begin{multline}
\label{eque}
G(x, u_{\epsilon}, Du_{\epsilon}, D^2u_{\epsilon})= \\
-\text{det}^{1/m} \left(\sigma^T\, D^2u\, \sigma+ Q(x, D u)+\nu (I_m+ \mu(x_1,\cdots, x_m)\otimes(x_1,\cdots, x_m))\right)+
\\ 
H^{1/m}(x, u_{\epsilon},  \sigma^T Du + \nu (x_1,\cdots, x_m)).
\end{multline}
From
 Minkowski's inequality \cite{HJ}:
\begin{equation}\label{Min}
\text{det}^{1/m}(A+B)\geq \text{det}^{1/m}(A)+\text{det}  ^{1/m}(B),\  A>0, B \geq 0, \text{ of order } m,
\end{equation}
\begin{multline
*}
G(x, u_{\epsilon}, Du_{\epsilon}, D^2u_{\epsilon})\leq 
-\text{det}^{1/m} (\sigma^T\, D^2u\, \sigma +Q(x, D u))-\\
\nu\,\text{det}^{1/m}(I_m+ \mu(x_1,\cdots, x_m)\otimes(x_1,\cdots, x_m))+
H^{1/m}(x, u_{\epsilon},  \sigma^T Du +\nu  (x_1,\cdots, x_m)).
\end{multline
*}
By Proposition \ref{stimagradiente}, $|\sigma^T(x)Du|\leq C$ in $\Omega_1$, so we use 
 the Lipschitz continuity of $H^{1/m}$ with $L=L_C$ and its monotonicity in $u$ \eqref{LipCH},
with the fact that $u$ is a subsolution of (\ref{MAsube0}), to obtain
$$
G(x, u_{\epsilon}, Du_{\epsilon}, D^2u_{\epsilon})\leq
\nu\left(L |(x_1,\cdots, x_m)|- \text{det} ^{1/m}\left(I_m+ \mu (x_1,\cdots, x_m)\otimes(x_1,\cdots, x_m)\right)\right).
$$
We want to find a suitably large $\mu$ independent of $\epsilon$  such that
\begin{equation}\label{stimaC}
 L^m |(x_1,\cdots, x_m)|^m< \det \big(I_m+ \mu (x_1,\cdots, x_m)\otimes(x_1,\cdots, x_m)\big), \mbox{ for any } x\in\Omega_1.
\end{equation}
We use the following equality \cite{ND}
\begin{equation}\label{NoDa}
\det(I+q\otimes q)= 1+|q|^2, \ q \mbox{ is a } m\times 1 \mbox{  column vector}.
\end{equation}
Then (\ref{stimaC}) becomes
\begin{equation}\label{ultima}
1+\mu |(x_1,\cdots, x_m)|^2 -  L^m |(x_1,\cdots, x_m)|^m>0.
\end{equation}
 \noindent i) If  $ |(x_1,\cdots, x_m)|<\frac{1}{L}$, then (\ref{ultima}) is true for any $\mu>0$;\\
  ii) if  $ |(x_1,\cdots, x_m)|\geq\frac{1}{L}$,  we take
  $$\mu >\max_{x}(L^m  |(x_1,\cdots, x_m)|^m -1)L^2,$$ 
  and 
  (\ref{ultima}) holds also in this case. 
  With this choice of $\mu$, for some $\a>0$ we have $G(x, u_{\epsilon}, Du_{\epsilon}, D^2u_{\epsilon})<-\a$
  for any $x\in\Omega_1$. Then \eqref{strict1} holds and the proof is complete.
  \end{proof}
  \begin{remark}\label{sigmaK}{\rm
  We can obtain the same result  also in 
  the case
 $$ \sigma(x)= \left(\begin{array}{cc}
K\\
\tau(x)\end{array}\right),$$
where $K$ is a nonsingular constant $m\times m$ matrix.
In this case 
 we use a generalization of (\ref{NoDa}):
 \begin{equation}\label{luigi}
\det \left(K^T(I+\mu v v^T)K\right)= \det (K^TK)\left(1+\mu |v|^2\right),\ 
\text{ if } |K^T v|\neq 0.
\end{equation}
 }
  \end{remark}
  
 \begin{remark}\label{sigmageneral} {\rm
 Another case where it is possible to costruct an uniformly $\chi$-convex strict subsolution is when
\begin{equation}\label{xsquare}
 \sigma^T(x)\sigma(x) +Q(x,x)\geq \eta\, I,\ \forall x\in \ol\Omega, \mbox{ for some } \eta>0.
 \end{equation}
  This condition is equivalent to say that $|x|^2$ is uniformly ${\mathcal X}$-convex in $\Omega$.
 In this case 
 we consider, for a given viscosity subsolution $u$ 
$$
u_{\epsilon}(x):=u(x)+\epsilon (e^{\mu\frac{|x|^2}{2}}-\lambda)
$$
and  we show that it is a strict subsolution 
for $\lambda >>\mu>>1$, independent of $\epsilon >0$,
 following  the procedure used in the proof of Theorem \ref{strict0}.
 To do this we apply the inequality \cite{PLL}:
\begin{equation}\label{disugPLL}
\det(A+\mu q\otimes q)\geq \eta^N(1+\frac{\mu}{N\eta}|q|^2),
\end{equation}
where
$A\in S^N$ such that $A\geq \eta\,I$, $\eta>0$ , $\mu\geq 0$, $q\in \R^N$.
}
\end{remark}


\subsection{The Comparison Principle with non-strict subsolutions}

  We are now ready to prove the main Comparison Principle for the equation
  \begin{equation}\label{E}
  -\det(D_ {\mathcal X}^2u)+H(x,u,D_{\mathcal X}u)=0,\ \mbox{in }\Omega.
  \end{equation}
\begin{theorem}\label{confrontogeneraleCarnot}
Assume $H$ satisfies \eqref{LipCH} and  
 the vector fields $X_1,...,X_m\in C^2$ are 
of Carnot type or satisfy \eqref{xsquare}.
 Let $u\in USC(\ol\Omega)$ be a bounded ${\mathcal X}$-convex subsolution of \eqref{E} 
and $v\in LSC(\ol\Omega)$ be a bounded supersolution of \eqref{E}. Then
  $$\sup_{\Omega}(u-v)\leq \max_{\partial \Omega}(u-v)^+.$$
  \end{theorem}
  \begin{proof}
We fix $\eta>0$ and set $s:= \max_{\partial \Omega}(u-v)^+$. By the upper semicontinuity of $u-v$ there is $\delta>0$ such that 
$$
(u-v)(x)\leq s+ \eta \quad\text{for all $x$ such that  dist}(x,\de\O)\leq\d.
$$
Now we set $\O_\d:=\{x\in\O\,:\, \text{dist}(x,\de\O)>\d\}$ and we must prove that $\sup_{\Omega_\d}(u-v)\leq s+ \eta$. To this goal we consider the uniformly ${\mathcal X}$-convex strict subsolutions $u_\e$ constructed in Theorem \ref{strict0}. We claim that each $u_\e$ is also strict subsolution of 
$$
- \log \det(D_ {\mathcal X}^2u_\e)+\log\left( H(x,u_\e,D_{\mathcal X}u_\e)+\frac\a2\right)\leq 0\; \mbox{ in }\Omega_\d ,
$$
where $\a>0$ is the constant provided by Theorem \ref{strict0}. Since $v$ is a supersolution of the same equation, we can then use Theorem \ref{confrontoC} in $\O_\d$ with $K_1=\log(H+\a/2)$ to get
$$
\sup_{\Omega_\d}(u_\e-v)\leq \max_{\partial \Omega_\d}(u_\e-v)^+\leq s+ \eta,\quad \forall\,\e\leq\e_o,
$$
where the last inequality follows from $u_\e\leq u$. Since $u_{\epsilon}\to u$, we let $\e\to 0$ and obtain that $u-v\leq s+ \eta$ in all $\O$, which gives the conclusion by the arbitrariness of $\eta$.

To prove the claim we recall that for each $\e$ there is $C$ such that $|D_{\mathcal X}u_\e|\leq C$ 
 in $\O_\d$, by Proposition \ref{stimagradiente}. Then 
 $$H(x,u_\e,D_{\mathcal X}u_\e)\leq \max_{x\in\ol\O, |p|\leq C}H(x, \sup u, p)=:C_1.
 $$
From this we get the  conclusion as in Lemma \ref{strictlogG}, because
\begin{multline*}
\log \left(H(x, u_\e, D_ {\mathcal X} u_\e) + \a\right) =
\log \left(H(x, u_\e, D_ {\mathcal X} u_\e)+\frac\a2\right)
+\log\left(1+\frac{\alpha}{2H(x, u_\e, D_ {\mathcal X} u_\e)+\a}\right)
\\
\geq \log \left(H(x, u_\e, D_ {\mathcal X} u_\e)+\frac\a2\right)+\log\left(1+\frac{\alpha}{2C_1+\a}\right).
\end{multline*}
 \end{proof}
 \begin{remark}
\label{no+2}
{\rm
If $H$ is independent of $u$ then the conclusion of the last theorem can be strengthened to
$\sup_{\Omega}(u-v)\leq \max_{\partial \Omega}(u-v)$
by Remark \ref{no+}.
}\end{remark}
 \begin{example}
 {\rm
 The assumption of the last theorem cover equations of the form
$$
-\det(D_{\mathcal X}^2u)+k(x,u)(1+|D_{\mathcal X}u|^2)^{\a}
=0,\ \mbox{ in } \Omega,
 $$
for  any $\a\geq 0$, $k\in C(\ol\O\times \R)$, $k\geq 0$ and nondecreasing in the second entry. If the vector fields are the canonical basis of the Euclidean space $\R^n$ and $\a= (n+2)/2$ this is the classical equation satisfied by a function $u$ whose graph has Gauss curvature $k$. In Carnot groups and for $\a= (m+2)/2$, $k=k(x)$, it is the equation of prescribed horizontal Gauss curvature as defined by Danielli, Garofalo, and Nhieu \cite{DGN}.
As a corollary of last theorem we obtain the uniqueness of a viscosity solution $u\in C(\ol\O)$ of this PDEs with prescribed boundary data.
}
\end{example}
\begin{example}
 {\rm
In problems of optimal transportation, see \cite{Ca, Vbook1} for the Euclidean case and  \cite{AR, FiRi},
$H$ has the form $H(x,q)=f(x)/h(q)$ with $f, h\geq 0$ and $\int_{\O}f(x)\,dx=\int_{\R^m} h(q)\, dq < +\infty$. The assumption \eqref{LipCH} of Theorem \ref{confrontogeneraleCarnot} is satisfied if $f\in C(\ol\O)$, 
$h\in C(\R^m)$, $h>0$, and $h^{-1/m}$ is locally Lipschitz. This is ok, for instance, if $h(q)=1/|q|^\a$ with $\a > m$.
}
\end{example}

In \cite{RT} Rauch and Taylor proved the following Comparison Principle for the classical Monge-Amp\`ere equation: \\
{\it
if $\O$ is strictly convex, $u\in C(\ol\O)$ is convex, $v\in W^{2,n}(\O)$, and $\det D^2u\geq\det D^2v$ in $\O$, then $\max_{\Omega}(u-v)\leq \max_{\partial \Omega}(u-v).$
}\\
The last result of this section, that is a special case of Theorem
\ref{confrontogeneraleCarnot} with Remark \ref{no+2}, gives a version of such statement in the context of viscosity solutions and noncommutative vector fields. It extends also a proposition of Gutierrez and Montanari \cite{GutMon} for $u,v\in C^2(\O)$, $n=3$, $m=2$, and $X_1, X_2$ generators of the Heisenberg group.

\begin{corollary}\label{corRT}
Assume the vector fields $X_1,...,X_m\in C^2$ are 
of Carnot type or satisfy \eqref{xsquare},
  $u\in USC(\ol\Omega)$  bounded and ${\mathcal X}$-convex, 
 $v\in LSC(\ol\Omega)$  bounded, and  
 $$
 -\det (D_{\mathcal X}^2u)+f(x)\leq 0, \quad -\det (D_{\mathcal X}^2v)+f(x)\geq 0 \quad\text{in } \O
 $$
 in viscosity sense for some $f\in C(\ol\O)$, $f\geq 0$.
  Then
  $\sup_{\Omega}(u-v)\leq \max_{\partial \Omega}(u-v).$
  \end{corollary}


 \section{Solvability of the Dirichlet problem}
 \label{DP}
 
In this section we apply the results of Section \ref{CPgeneral}
to solve the 
Dirichlet problem for the PDE \eqref{E}
\begin{equation}\label{DP1}
\left\{\begin{array}{ll}
-\det 
D_{\mathcal X}^2 u
+ H(x, u, D_{\mathcal X}u)=0\quad &\mbox{ in } \ \Omega ,\\
u=g\  & \mbox{ on }\ \partial\Omega,
\end{array} \right.
\end{equation}
with $g\in C(\de\Omega)$.
\subsection{Some explicit solutions in the Heisenberg group}
\label{explicit}
In $\R^3$ with coordinates $x=(x_1,x_2,t)$ the generators of the Heisenberg group are the vector fields
\begin{equation}
\label{genH}
X_1=\frac{\partial}{\partial x_1} + 2x_2\frac{\partial}{\partial t}, \qquad X_2=\frac{\partial}{\partial x_2} - 2x_1\frac{\partial}{\partial t} .
\end{equation}
The norm 
\begin{equation}
\label{normH}
|x|_\mathbb{H}:= w(x)^{1/4}, \qquad w(x_1,x_2,t):= (x_1^2 + x_2^2)^2 + t^2 ,
\end{equation}
is positively 1-homogeneous with respect to the dilations $\d_\lambda(x_1,x_2,t)=(\lambda x_1,\lambda x_2,\lambda^2t)$.
The Koranyi ball of radius $R>0$ centered at the origin is 
\[
B_\mathbb{H}(R):=\left\{x=(x_1,x_2,t)\in\mathbb{R}^3 \,:\, |x|
_\mathbb{H}
<R
\right\}.
\]
\begin{proposition}
\label{explicit1}
Let $\O = B_\mathbb{H}(R)$ and $X_1,X_2$ be the generators of the Heisenberg group \eqref{genH}.
Then $w(x)=|x|^4_\mathbb{H}$ 
 is the unique $\mX$-convex viscosity solution of the Dirichlet problems
\begin{equation}
\label{MAHeis}
-\det D_{\mathcal X}^2 u + 144\left(x_1^2 + x_2^2\right)^2 = 0 \;\;\mbox{ in }  \Omega ,\quad
u=R^4\;  \mbox{ on }\; \partial\Omega,
\end{equation}
and 
\begin{equation}
\label{GaussHeis}
-\det D_{\mathcal X}^2 u + k_\mathbb{H}(x)\left(1+|D_{\mathcal X}u|^2\right)^2  = 0 \;\;\mbox{ in }  \Omega ,\quad
u=R^4\;  \mbox{ on }\; \partial\Omega,
\end{equation}
where
\begin{equation}
\label{curvHeis}
 k_\mathbb{H}(x):=\left(\frac{12\left(x_1^2 + x_2^2\right)
 }{
 1+16(x_1^2 + x_2^2)|x|^4_\mathbb{H}}\right)^2 .
\end{equation}
 In particular, $w$ is the unique $\mX$-convex function on $B_\mathbb{H}(R)$ with 
 horizontal Gauss curvature  $k_\mathbb{H}$ and boundary value $R^4$.
\end{proposition}
\begin{proof}
The uniqueness follows from Theorem \ref{confrontogeneraleCarnot}.
A straightforward calculation gives
\begin{equation}
\label{hessw}
 D_{\mathcal X}^2 w(x) = 12\left(x_1^2 + x_2^2\right) I ,
\end{equation}
so $w$ is $\mX$-convex and a classical solution of \eqref{MAHeis}.
Moreover 
\begin{equation}
\label{gradw}
D_{\mathcal X} w = 4\left(x_1^2 + x_2^2\right)
 \left(\begin{array}{c}x_1\\ x_2\end{array}\right) 
 + 4t \left(\begin{array}{c}x_2\\ -x_1\end{array}\right) ,
\end{equation}
so 
\begin{equation}
\label{|gradw|}
|D_{\mathcal X} w|^2 = 16(x_1^2 + x_2^2)\left((x_1^2 + x_2^2)^2 + t^2\right)
\end{equation}
and $w$ is  a classical solution of \eqref{GaussHeis}.
\end{proof}

The next result is the analogue in the Heisenberg group of the fact that the Euclidean norm $|x|$ in $\R^n$ is the unique convex function solving $\det D^2 u=0$ in the punctured Euclidean ball $B(R)\setminus \{0\}$ and taking the values $R$ on $\partial B(R)$ and $0$ at $0$. 
\begin{proposition}
\label{explicit2}
Let 
$X_1,X_2$ be the generators of the Heisenberg group \eqref{genH}.
Then the homogeneous norm $|\cdot|_\mathbb{H}$ 
 is the unique $\mX$-convex viscosity solution of the Dirichlet problem
\begin{equation}
\label{MAHeis2}
-\det D_{\mathcal X}^2 u  = 0 \quad\mbox{in }  B_\mathbb{H}(R) \setminus \{0\},\quad
u=R \quad 
 \mbox{on }\partial B_\mathbb{H}(R),\quad u(0)=0.
\end{equation}
\end{proposition}
\begin{proof}
The uniqueness follows from Theorem \ref{confrontogeneraleCarnot}. Next we compute, for $x\ne 0$,
\begin{equation}
\label{D2norm}
D^2_{\mathcal X}|x|_\mathbb{H}  = \frac{1}{4|x|_\mathbb{H}^3}\left[D^2_{\mathcal X}w -  \frac{3}{4w}D_{\mathcal X}w\otimes D_{\mathcal X}w\right].
\end{equation}
If $(x_1, x_2)=0$, then $D^2_{\mathcal X}|x|_\mathbb{H}=0$. 
If  $(x_1, x_2)\neq 0$, 
to show that the matrix in brackets $[...]$ is positive semidefinite we take a unit vector $\zeta$ and set $\psi:=x_1^2 + x_2^2.$ First observe that
\begin{equation*}
\zeta^T D^2_{\mathcal X}w \, \zeta = 12 \psi |\zeta|^2 = 12 \psi
\end{equation*}
by \eqref{hessw}. Next we compute, using \eqref{|gradw|},
\begin{equation*}
\frac{3}{4w}\zeta^T \left(D_{\mathcal X}w\otimes D_{\mathcal X}w\right)\, \zeta = \frac{3}{4w}\left|\zeta \cdot D_{\mathcal X}w\right|^2 \leq \frac{3}{4w}\left|D_{\mathcal X}w\right|^2 = \frac{3}{4w} 16 \psi w = 12 \psi .
\end{equation*}
Then $D^2_{\mathcal X}|x|_\mathbb{H}\geq 0$ in $\R^3\setminus\{0\}$ in the classical sense and $|\cdot|_\mathbb{H}$  is  $\mX$-convex.

To prove that $\det D^2_{\mathcal X}|x|_\mathbb{H}=0$ it is enough to show, by \eqref{D2norm}, that
\begin{equation}
\label{=0}
\left[D^2_{\mathcal X}w -  \frac{3}{4w}D_{\mathcal X}w\otimes D_{\mathcal X}w\right] D_{\mathcal X} w = 0,
\end{equation}
because $D_{\mathcal X} w\ne 0$ for $(x_1, x_2)\neq 0$. By \eqref{hessw} and  \eqref{gradw} the first term is
\begin{equation*}
D^2_{\mathcal X}w  D_{\mathcal X} w = 48 \psi^2 
 \left(\begin{array}{c}x_1\\ x_2\end{array}\right) 
 + 48 t \psi \left(\begin{array}{c}x_2\\ -x_1\end{array}\right) .
\end{equation*}
For the second term we use $\left(D_{\mathcal X}w\otimes D_{\mathcal X}w\right) D_{\mathcal X} w = \left| D_{\mathcal X} w \right|^2 D_{\mathcal X} w$, \eqref{gradw}, and finally \eqref{|gradw|} to compute
\[
\frac{3}{4w}\left(D_{\mathcal X}w\otimes D_{\mathcal X}w\right) D_{\mathcal X} w = \frac{3}{4w} 16\psi w D_{\mathcal X} w = 48 \psi^2 
 \left(\begin{array}{c}x_1\\ x_2\end{array}\right) 
 + 48 t \psi \left(\begin{array}{c}x_2\\ -x_1\end{array}\right) ,
\]
which gives \eqref{=0}.
\end{proof}
\begin{remark}
{\rm
The calculations of this section hold as well in $\R^{2j+1}$ with the generators of the $j$-th Heisenberg group and the corresponding homogeneous norm. The fact that such norm solves $\det D^2_{\mathcal X}=0$ off the origin was proved in \cite{DGN} for general stratified groups of Heisenberg type.
}\end{remark}
\subsection{Perron method}
\label{per}
Here we describe the construction of solutions in the general case.
We denote by $\mathcal S$ and $\mathcal Z$, respectively, the sets of 
sub- 
and supersolutions 
of $(\ref{DP1})$: 
\begin{eqnarray*}
&\mathcal S:=\{ w\in USC(\ol \Omega): \mbox {$w$  bounded $\mathcal X$-convex subsolution of } \eqref{E}
,\ w\leq g, \mbox { on } \partial\Omega
\},\\
&\mathcal Z:=\{ W\in LSC(\ol \Omega): \mbox {$W$  bounded supersolution of } \eqref{E} 
,\ W\geq g,\mbox { on } \partial\Omega
\}.
\end{eqnarray*}
The Perron 
method proposes as a candidate solution of \eqref{DP1}
$$
\ul u(x):=\sup_{w\in \mathcal S}w(x),\ x\in \ol\Omega,\quad\text{if } \mathcal S \ne \emptyset .
$$
Note that, if $W\in \mathcal Z$ and the Comparison Principle holds, then
$\ul u(x)\leq W(x) < +\infty$ for all $x$. Under no further assumptions $\ul u$ is a generalized, possibly discontinuous solution of the Dirichlet problem \eqref{DP1} in the following sense.
\begin{theorem}
\label{PWB}
If $\mathcal S\ne \emptyset$, $\mathcal Z\ne \emptyset$, and the Comparison Principle holds for \eqref{DP1}, then the u.s.c. envelope $\ul u^*$ is a $\mathcal X$-convex subsolution and the l.s.c. envelope $\ul u_*$ is a supersolution of \eqref{E}.
\end{theorem}
 \begin{proof}
  By a standard argument, if a function $v$ is the sup of a set of subsolutions then its u.s.c. envelope $v^*(x):=\inf\{V(x)\, :\, V\in USC(\ol\Omega), v\leq V\}$ is a subsolution, see \cite{CIL, BCD}. In particular,
 $\ul u^*(x)$ is $\mathcal X$-convex. 
 
 The proof that $\ul u_*$ is a supersolution of \eqref{E} is achieved by contradiction: one assumes that $\ul u_*$ fails to be a supersolution at some point $y\in\Omega$ and constructs a subsolution that is larger than $\ul u$ near $y$, therefore contradicting the maximality of $\ul u$. For the Monge-Amp\`ere equations we must show that the we can construct a $\mathcal X$-convex subsolution larger than $\ul u$.

 If $\ul u_*$ fails to be a supersolution at some point, say $y=0$, by Definition \ref{defsupersol} there is a $C^2$ test function $\varphi$ such that $\varphi(0)=\ul u_*(0)$, $\varphi(x)\leq \ul u_*(x)$ for $|x|$ small,
$$
 -\det (D_{\mathcal X}^2 \varphi(0))+ H(0, \ul u_*(0), D_{\mathcal X}\varphi (0))<0 ,
 $$
 and 
 \begin{equation}
 \label{strictphi}
 D_{\mathcal X}^2 \varphi(0)>0.
 \end{equation}
 The usual "bump" construction \cite{CIL, BCD} considers $v(x):=\varphi(x)+\e-\gamma |x|^2$ and 
 $$
 U(x):=\max\{\ul u(x), v(x)\}\;\text{ if } |x|<r, \quad U(x):=\ul u(x)\;\text{ otherwise.}
 $$
 Then one checks that for small $\e, \g, r$, $U^*$ is a subsolution of the PDE \eqref{E} and $\sup(U^*-\ul u)>0$. Thanks to \eqref{strictphi}, by further restricting $\g, r$ if necessary, we also have that $-D_{\mathcal X}^2 U^*\leq 0$ in viscosity sense, so $U^*$ is $\mathcal X$-convex and we achieve the contradiction by the usual argument \cite{CIL, BCD}.
   \end{proof}
  
  To give examples of equations that satisfy the last theorem, as well as the next results on the existence of continuous solutions, we will use the assumption that for some $L, M, R>0$
 \begin{equation}
 \label{crescita}
H^{1/m}(x,R,p)\leq L|p| + M \quad \forall\; x\in\ol\O,\; p\in\R^m .
\end{equation}
Since $H$ is nondecreasing in the second entry, the same inequalitiy holds for $H^{1/m}(x,r,p)$ and all $r\leq R$. In the next results we will 
take $R=\min_{\de\Omega}g$ or $\max_{\de\Omega}g$, and $L, M$ will depend on it.
This is a slightly weaker version of the growth condition used by P.-L. Lions \cite{PLL2} in the Euclidean case, see Examples \ref{crescitapll} and \ref{megliopll} below.
   \begin{example}
  \label{soluzgeneralizCarnot}
  {\rm
Assume $H$ satisfies \eqref{LipCH}, 
  \eqref{crescita} with $R=\min_{\de\Omega}g$,  
and the vector fields $X_1,...,X_m\in C^2$ are 
of Carnot type. 
Then $\ul u(x)$ is a generalized solution of problem (\ref{DP1}) in the sense described by the last theorem. 
%
In fact, by Theorem \ref{confrontogeneraleCarnot}
we know that the Comparison Principle 
for (\ref{DP1}) holds.
Therefore, it is enough to prove that both sets 
$\mathcal S$ and  $\mathcal Z$ are nonempty.
First of all we note that $W\equiv \max_{\de\Omega}g$ is 
an element of $\mathcal Z$.
As far as the set $\mathcal S$, we consider
 $$
 w(x)= e^{\frac{\mu}{2}\sum_{i=1}^m |x_i|^2}-\max_{x\in\ol\Omega}e^{\frac{\mu}{2}\sum_{i=1}^m |x_i|^2} + \min_{\de\Omega}g ,
 $$
so that $w\leq\min_{\de\Omega}g$.
As in the proof of Theorem \ref{strict0},  we have that  $D^2_{\mathcal X}w\geq \mu I$. Moreover, for $\nu:= \mu e^{\mu\frac{\sum_{i=1}^m |x_i|^2}{2}}$, by \eqref{NoDa}  we get
\begin{multline*}
-\text{det}^{1/m} (D^2_\mX w)+
H^{1/m}(x, w,  D_{\mX}w)=\\
-\nu \text{det}^{1/m} \left(I_m+ \mu(x_1,\cdots, x_m)\otimes(x_1,\cdots, x_m)\right)+
H^{1/m}(x, w,  \nu (x_1,\cdots, x_m)
)\\
\leq
\nu \left(-1-\mu|(x_1,\cdots, x_m)|^{2/m} + L|(x_1,\cdots, x_m)|\right) +M
\nn
\end{multline*}
which becomes negative for $\mu$ large enough.
}\end{example}
  \begin{example}
  \label{soluzgenerxsquare}
  {\rm
Assume $H$ satisfies \eqref{LipCH}, 
  \eqref{crescita} with $R=\min_{\de\Omega}g$, and the vector fields are of class $C^2$ and such that $|x|^2$ is ${\mathcal X}$-convex in $\Omega$, that is, the inequality \eqref{xsquare} holds. 
Then  $\ul u(x)$ is a generalized solution of problem (\ref{DP1}) as in the preceding example. The proof is the same except that now we use
$
 w(x)=e^{\mu\frac{|x|^2}{2}}- \max_{x\in\ol\Omega}e^{\mu\frac{|x|^2}{2}}+ \min_{\de\Omega}g .
 $
}\end{example}

As in the classical potential theory, the continuity at the boundary of the Perron solution requires the existence of barriers.
\begin{definition}
We say that $w$ is a lower (respectively, upper) barrier for problem (\ref{DP1}) at a point 
$x\in\partial\Omega$
if $w\in \mathcal S$ (respectively, $w\in \mathcal Z$) and
$$\lim_{y\rightarrow x}w(y)= g(x).$$
\end{definition}
\begin{corollary
}\label{PWB2} 
Suppose that the Comparison Principle holds for 
\eqref{DP1} and that for all $x\in \partial\Omega$ there exist a lower and an upper barrier.
Then 
$\ul u\in C( \ol\Omega)$ is the solution of \eqref{DP1}, that is, the unique $\mX$-convex viscosity solution of \eqref{E} attaining continuously the boundary data $g$.
\end{corollary
}
\begin{proof}
The existence of a lower and an upper barrier at $x\in\de\Omega$ implies the continuity of $\ul u$ at $x$ and $\ul u(x)=\ul u^*(x)=\ul u_*(x)=g(x)$. Then the Comparison Principle gives $\ul u^*=\ul u_*$ in $\Omega$ and therefore $\ul u$ is a continuous viscosity solution of \eqref{E}.
\end{proof}
\begin{remark}
{\rm {\em Interior regularity of the solution.}
Since the solution $\ul u$ of the Dirichlet problem is $\mX$-convex, it is locally Lipschitz continuous with respect to the Carnot-Carth\'eodory distance $d$ associated to the vector fields $\mX$ (see, e.g., \cite{Bell, BD} for the definition). 
If, in addition, the identity map $(\R^n,d) \to  (\R^n,|\cdot|)$ is a homeomorphism (e.g., the vector fields $\mX$ are smooth and satisfy the H\"ormander condition), then the distributional derivatives $X_j\ul u$ exist a.e. and are locally bounded. All this is known in Carnot groups by \cite{LMS, BaRic,DGN, Ric, Mag} and was proved in \cite{BD} for general vector fields. 

In Carnot groups of step 2, horizontally convex functions are also twice differentiable a.e. \cite{GutMon2, DGNT, Mag}. Therefore in this case the Perron solution $\ul u$ solves the PDE \eqref{E} also pointwise almost everywhere.
}\end{remark}

%
\subsection{Construction of barriers}
\label{bar}

To find explicit examples where the Perron method works and the Dirichlet problem is solvable we make some assumptions on the bounded open set $\Omega$. We say it is smooth if
\begin{eqnarray}
&&\mbox{there exists 
}\ \Phi
\in C^2\ \mbox{such that}\label{rappr}\\
&&\Omega=\{x\in \re^n: \Phi(x)>0\},\ 
D\Phi(x)\ne0,\ \forall x\in \partial \Omega.\nn
\end{eqnarray}
%
The main additional  assumptions is that the domain be uniformly 
convex with respect to the vector fields $X_j$. 
It is the natural extension for subelliptic Monge-Amp\`ere equations of the standard uniform convexity 
in the Euclidean case \cite{GT, PLL, PLL2, CNS, Tru06}. 
\begin{definition}
A domain $\Omega$ smooth in the sense of (\ref{rappr}) is called 
 convex with respect to the fields 
  $X_1,\cdots,X_m$ (briefly,  $\mathcal X$-convex) and, respectively, uniformly $\mathcal X$-convex, if 
 $D^2_{\mX}\Phi(x)\leq 0$ and, respectively, 
\begin{equation}\label{domainconvex}
D^2_{\mX}\Phi(x)\leq -\gamma I
\quad\text{for some $\gamma>0,$} \mbox{ for any }\ x\in\Omega.
\end{equation}
\end{definition}

\begin{example}{\rm
Any Euclidean ball centered in $x_o$ is uniformly $\mathcal X$-convex if and only if 
$|x-x_o|^2$ is ${\mathcal X}$-convex in $\R^n$.  It is well known that this is true in all Carnot groups of step 2, in particular the Heisenberg groups.
 For $x_o=0$ the $\mathcal X$-convexity of $|x|^2$ is equivalent to the inequality \eqref{xsquare} 
in $\R^n$ that we already used in the Comparison Principles.

The Koranyi ball $B_\mathbb{H}(R)$ in $\R^3$ is $\mX$-convex but not uniformly $\mX$-convex with respect to the generators of the Heisenberg group \eqref{genH}.
}
\end{example}

The next result is an analogue in the context of Carnot-type 
vector fields of a classical result of Caffarelli, Nirenberg, and Spruck \cite{CNS} and P.-L. Lions \cite{PLL2} in the Euclidean case, saying that the solvability of the Dirichlet problem is equivalent to the existence of a convex subsolution attaining continuously the boundary data (see also \cite{Tru90} for other curvature equations).

\begin{corollary}\label{cnsl} Assume $H$ satisfies (\ref{LipCH}) and 
the vector fields $X_1,...,X_m\in C^2$ are either 
of Carnot type or they 
satisfy \eqref{xsquare}. Suppose also that either (i) $g\equiv 0$, or (ii) $\O$ and $g\in C^2(\ol\Omega)$ are 
$\mX$-convex, or (iii) $\O$ is uniformly $\mX$-convex and $g\in C^2(\ol\Omega)$. If for all $x\in \partial\Omega$ there exist a lower barrier, then $\ul u\in C( \ol\Omega)$ is the unique solution of \eqref{DP1}.
\end{corollary}
\begin{proof}
We are going to apply  Corollary \ref{PWB2}. The Comparison Principle comes from  Theorem \ref{confrontogeneraleCarnot}. An upper barrier at all points of the boundary is $W\equiv 0$ in
case $(i)$. In the other cases we take $W(x):= \lambda \Phi(x) - g(x)$, $\lambda >0$. Then in case $(ii)$ $D^2_{\mX}W(x)\leq 0$ for any $\lambda$ and in case $(iii)$ $D^2_{\mX}W(x)\leq - \lambda \gamma I - D^2_{\mX}g(x)\leq 0$ 
for  $\lambda$ large enough. Therefore any test function $\varphi$ such that $W-\varphi$ attains a minimum at $x$ has   $D^2_{\mX}\varphi(x)\leq 0$. Then $W$ is a supersolution of \eqref{E} by Definition \ref{defsupersol}.
\end{proof}

\begin{remark}{\rm
In the papers  \cite{CNS} and \cite{PLL2}  the vector fields are the canonical basis of $\R^n$ and it is assumed the existence of a convex subsolution of class  $C^2$ (or solving the PDE in the sense of Alexandrov) that is also a  lower barrier at all points. On the other hand it is proved the existence of a $C^\infty$ solution to the Dirichlet problem.
}
\end{remark}

Next we give some explicit conditions on the data ensuring the existence of a lower barrier and therefore the solvability of the Dirichlet problem.

\begin{proposition}\label{barrierainfglobale}
Suppose that $\Omega$ is 
uniformly $\mathcal X$-convex and $g\in C^2(\ol\Omega)$.
Assume $H$ satisfies
 \eqref{crescita} with $R=\max_{\partial \Omega}g$ and it is nondecreasing w.r.t. the second entry $r$.
Then 
there exists 
$w\in {\mathcal S}\cap C(\ol\O)$ such that $w=g$ on $\partial\O$.
\end{proposition}
\begin{proof}
We consider
\begin{equation}\label{barrieraglo}
w(x)=\lambda(e^{-\mu \Phi(x)}-1) +g(x),\quad \mu, \lambda >0,
\end{equation}
where $\Phi$ is defined in (\ref{rappr}). 
Clearly $w(x)=g(x)$ for any $x\in\partial\Omega$ and $w(x)<g(x)$ for any $x\in\Omega$. Moreover
\begin{eqnarray*}
&&Dw(x)=-\mu\lambda\,e^{-\mu \Phi}D\Phi+Dg ,\nn\\
&&D^2w(x)=\mu\lambda\,e^{-\mu \Phi}
(-D^2\Phi +\mu D\Phi \otimes D\Phi) +D^2g ,\nn\\
&&D_{\mathcal X}^2 w= \mu\lambda\,e^{-\mu \Phi}
(-D_{\mathcal X}^2\Phi +\mu \sigma^TD\Phi \otimes \sigma^TD\Phi) +D_{\mathcal X}^2g .
\end{eqnarray*}
From the regularity of the function $g$
(there is $c$ such that $-c\,I\leq D_{\mathcal X}^2g$) 
and the uniform $\mathcal X$-convexity of $\Omega$ with constant $\gamma$
we have
$$D_{\mathcal X}^2 w\geq  (\mu\lambda\,e^{-\mu \Phi}\gamma-c)I
+\mu^2\lambda\,e^{-\mu \Phi}\sigma^TD\Phi \otimes \sigma^TD\Phi. $$
First we choose $\lambda$ such that
$
\mu\lambda\,e^{-\mu \Phi(x)}\gamma-c \geq  \frac{1}{2}\mu\lambda\,e^{-\mu \Phi(x)}\gamma,
$ 
i.e.
 \begin{equation}\label{primolambda}
\lambda \geq  \frac{2c}{\mu
\gamma}\max_{\ol\O}e^{\mu \Phi}.
\end{equation}
With this $\lambda$, $D_{\mathcal X}^2 w\geq 0$
and $w$ is $\mathcal X$-convex.
Moreover, from inequality (\ref{disugPLL}),
\begin{eqnarray}\label{detbarr}
&&\det(D_{\mathcal X}^2 w)\geq (\mu\lambda\,e^{-\mu \Phi})^m \det (\frac{\gamma}{2}\, I
+\mu\sigma^TD\Phi \otimes \sigma^TD\Phi) \\
&&\geq
 (\mu\lambda\,e^{-\mu \Phi})^m \left(\frac{\gamma}{2}\right)^m
\left(1
+\frac{2\mu}{m\gamma}|\sigma^TD\Phi|^2\right).\nn
\end{eqnarray}
From the monotonicity of $H$, the growth assumption \eqref{crescita} with $R=\max_{\partial \Omega}g$, and the boundedness of $Dg$ we get, for some $K>0$,
\begin{eqnarray}\label{Hbarr}
&&H(x,w, D_{\mathcal X}w)\leq H(x,R, D_{\mathcal X}w)\leq (M+ L | D_{\mathcal X}w+ D_{\mathcal X}g |)^m\\
&&\leq 
K+ K (\mu\lambda\,e^{-\mu \Phi(x)})^m\, |\sigma^TD\Phi|^m.\nn
\end{eqnarray}
To prove that $w$ is a subsolution we have to show that the r.h.s. of \eqref{detbarr} is larger than the r.h.s. of \eqref{Hbarr}
for $\mu$ and $\lambda$ large enough.
This  is equivalent to
$$
1+\frac{2\mu}{m\gamma}|\sigma^TD\Phi|^2\geq
\frac{2^m K}{ (\gamma\mu\lambda\,e^{-\mu \Phi})^m} + \left(\frac{2}{\gamma}\right)^m K |\sigma^TD\Phi|^m.$$
First  we choose $\mu$ so large that
$$ 
\frac{2\mu}{m\gamma}|\sigma^TD\Phi|^2\geq \left(\frac{2}{\gamma}\right)^m K |\sigma^TD\Phi|^m
$$
and then we choose
$\lambda$ such that  $2^m K/ (\gamma\mu\lambda\,e^{-\mu \Phi})^m\leq 1$,
i.e.,
$$
\lambda \geq
 \frac{2K^{1/m}}{\mu\gamma}\max_{\ol\O}e^{\mu \Phi}.
$$
From (\ref{primolambda}) we can conclude by choosing
$$
\lambda \geq
\frac{2}{\mu\gamma}\max(c, K^{1/m})\max_{\ol\O}e^{\mu \Phi}.$$
\end{proof} 
%
\begin{theorem}\label{soluzcontinuaCarnot}
Suppose that  $\Omega$ is uniformly $\mathcal X$-convex, $g\in C(\partial \O)$,
and  the vector fields $X_1,...,X_m\in C^2$ are either of Carnot type or they   satisfy \eqref{xsquare}.
Assume $H$ satisfies \eqref{LipCH} and  \eqref{crescita} with $R=\max_{\partial \O} g$.
Then $\ul u\in C( \ol\Omega)$ is the unique solution of \eqref{DP1}.
 \end{theorem}
 \begin{proof}
 We take a sequence $g_n\in C^2(\ol \O)$ that converges uniformly to $g$. By Corollary \ref{cnsl} and Proposition \ref{barrierainfglobale} there is a solution $u_n\in C(\ol\O)$ of the Dirichlet problem with boundary condition  $g_n$.  By the estimate of Theorem \ref{confrontogeneraleCarnot} 
 $$
 \sup_{\Omega}|u_n-u_m|\leq \max_{\partial \Omega}|g_n-g_m| \quad \forall\, n, m. 
 $$
Since  $g_n$ is a Cauchy sequence also $u_n$ is such and therefore it converges uniformly to $u\in C( \ol\Omega)$. By the stability of viscosity solutions $u$ solves \eqref{DP1} and by the Comparison Principle Theorem \ref{confrontogeneraleCarnot} it coincides with $\ul u$.
\end{proof}

Next we list several equations to which the last theorem on the well-posedness of the Dirichlet problem applies.
\begin{example}
\label{crescitapll}
{\rm
Assume the Hamiltonian $H$ in \eqref{DP1} satisfies \eqref{LipCH} and for some $R\in \R$
$$
\limsup_{|p|\to +\infty} \frac{H(x,R,p)}{|p|^m} <+\infty\quad\text{uniformly in $x\in\ol \Omega.$}
$$
Then 
the growth assumption \eqref{crescita} holds for 
some $L, M$ 
and Theorem \ref{soluzcontinuaCarnot} applies for all data $g\in C(\partial \O)$ such that 
$\max_{\partial \O} g \leq R$.
In the Euclidean case ($m=n$ and the vector fields are the canonical basis of $\R^n$) we therefore recover one of the main results of P.-L. Lions' paper \cite{PLL2}. 
}
\end{example}
\begin{example}
{\rm
The previous example includes the basic subelliptic Monge-Amp\`ere equation \cite{Man, GutMon, DGN}
$$
-\det (D_{\mathcal X}^2 u)+f(x)=0,  \quad f\in C(\ol\Omega), \quad f\geq 0,
$$
as well as  the case 
 $H(x,r,p)=(1+\lambda |p|^2)^{m/2}$, with $\lambda>0$, that in the Euclidean case was studied in \cite{PLL}.
}
\end{example}
\begin{example}
{\rm
Equations of the form
$$
-\det (D_{\mathcal X}^2 u)+f(x, u)=0,  \quad f\in C(\ol\Omega\times \R) \text{ nondecreasing in }u, \quad f\geq0,  
$$
satisfy the assumptions of Theorem \ref{soluzcontinuaCarnot}. The dependence of $f$ on $u$ appears in various problems of differential geometry, see \cite{Au, CY77} and the references therein. An example that appears in several papers \cite{Del, CY82, Au} is
$$
f(x, u)=\phi(x)e^{\lambda(x) u}, 
 \quad \lambda, \phi\in C(\ol\Omega), \quad  \lambda, \phi \geq 0.
$$
}\end{example} 
\begin{example}
\label{megliopll}
{\rm
Another special case of Example \ref{crescitapll} is
\begin{equation*}
-\det (D_{\mathcal X}^2 u) + k(x,u) \left(1+|D_{\mathcal X} u|^2\right)^\alpha
 =0,\quad k \geq 0,\; 
 \; 0\leq\alpha\leq \frac m2, 
\end{equation*}
with $k\in C(\ol\Omega\times \R)$  nondecreasing in $u$, whose structure is reminiscent of the equation of prescribed Gauss curvature (where, however, $\alpha=\frac m2+1$). We recall that the exponent $\frac m2$ is optimal without additional compatibility conditions on $k$. 
In the Euclidean case it is well-known that a necessary condition for the existence of a classical subsolution attaining the boundary data is 
$$
\int_\Omega k(x)\, dx \leq \int_{\R^n} (1+|p|^2)^{-\alpha} \,dp ,
$$
and the right hand side is finite if and only if $\alpha> n/2$.
}
\end{example}
\begin{example}
{\rm
A simple example of nonexistence of the viscosity solution to the Dirichlet problem when the growth assumption \eqref{crescita} fails is the ODE
$$
-u''+H(x,u,u')=0 \quad\text{in }]-1,1[, \quad\text{with } H(x,u,p)\geq\frac\pi 2\left(1+p^2\right),
$$
with boundary conditions $u(-1)=u(1)=0$. Indeed there are no viscosity subsolution, i.e., $\mathcal S=\emptyset$. To prove this we solve
$-u''+k\left(1+(u')^2\right)=0$ with $k<\frac\pi 2$ and find $u_k(x)=\frac 1k \log\frac{\cos k}{\cos kx}$. If $w\in\mathcal S$ the Comparison Principle gives $w\leq u_k$, but $u_k(x)\to-\infty$ as $k\to \frac\pi 2 -$, for all $|x|<1$.
}
\end{example}

\begin{remark}\label{charac}
{\rm
In the theory of linear 
PDEs it is important to distinguish the characteristic and non-characteristic points of the boundary. For fully nonlinear, degenerate elliptic PDEs 
$$F(x,u,Du, D^2u)=0\quad \text{in } \Omega$$ we gave in \cite{BM} the following definition:

\noindent a point $z\in\partial\Omega$ is called {\emph{characteristic}} for the operator $F$ if 
\begin{eqnarray}
F(z, 0, -n(z), X+\mu n(z)\otimes n(z))=F(z, 0, -n(z), X),\
\forall X,\ \forall \mu >0,\phantom{t}\label{PCO} 
\end{eqnarray}
For the subelliptic Monge-Amp\`ere 
equations \eqref{E} the operator $F$ is
\begin{equation*}\label{Eq2}
F(x,r, p,X):= - \det  (\sigma^T X\sigma +Q(x,p))+ H(x,r, \sigma^T p),
\end{equation*}
with $Q$ given by (\ref{Q}).
The characteristic points $z$ are determined by the equation
\begin{equation}\label{PCeq}
\det (\sigma^T(z) X\sigma(z) +Q(z,p)+\mu \sigma^T(z)n(z)\otimes \sigma^T(z)n(z))=
\det (\sigma^T(z) X\sigma(z) +Q(z,p)).
\end{equation}
We recall that if $A$, $B$  are square matrices, $A>0$ and $B\geq 0$, then
$$
\det(A+B)=\det (A)\ \text{if and only if}\  B\equiv 0.$$
Then (\ref{PCeq}) for  $\mu>0$ yields  $ |\sigma^T(z)n(z)| =0$. In the proof of Theorem \ref{soluzcontinuaCarnot} we did not need to treat these points differently form the non-characteristic ones 
because of the uniform $\mX$-convexity of the domain. For more general domains we expect that the lower barrier must be constructed in different ways at characteristic and non-characteristic points, as for the linear 
and the Hamilton-Jacobi-Bellman subelliptic PDEs, see \cite{BLU, BM} and the references therein.
}
\end{remark}

We end the section with two results in the Koranyi ball of the Heisenberg group. Note that $B_\mathbb{H}(R)$ is  a $\mX$-convex set but it is not uniformly $\mX$-convex, because the condition \eqref{domainconvex} fails at the characteristic points of the boundary, namely, the intersections with the $t$ axis.
\begin{corollary}
Let 
 $X_1,X_2$ be the generators of the Heisenberg group \eqref{genH}.
Then, for any $f\in C(\ol B_\mathbb{H}(R)\times \R)$ nondecreasing in the second entry and such that  $f(x,R^4)\leq 144(x_1^2 + x_2^2)^2$, $\ul u\in C(\ol B_\mathbb{H}(R))$ 
 is the unique $\mX$-convex viscosity solution of the Dirichlet problem
\begin{equation}
\label{MAHeisf}
-\det D_{\mathcal X}^2 u + f(x,u) = 0 \;\;\mbox{ in }  B_\mathbb{H}(R)
 ,\quad
u=R^4\;  \mbox{ on }\; \partial B_\mathbb{H}(R).
\end{equation}
\end{corollary}
\begin{proof}
By Proposition \ref{explicit1} $w(x)=|x|^4_\mathbb{H}$ is a subsolution of \eqref{MAHeisf} that attains continuously the boundary data.
Then the conclusion follows from Corollary \ref{cnsl} .
\end{proof}
The last result is about the prescribed horizontal Gauss curvature equation. Although it does not satisfy the growth condition \eqref{crescita}, a lower barrier is given by $w(x)=|x|^4_\mathbb{H}$ 
if the prescribed curvature is lower than the horizontal curvature $k_\mathbb{H}$ of the graph of $w$.
\begin{corollary}
Let 
 $X_1,X_2$ be the generators of the Heisenberg group \eqref{genH}.
Assume $k\in C(\ol B_\mathbb{H}(R)\times \R)$ is nondecreasing in the second entry and satisfies  $k(x, R^4)\leq k_\mathbb{H}(x)$, where $k_\mathbb{H}$ is given by \eqref{curvHeis}.
Then $\ul u\in C(\ol B_\mathbb{H}(R))$ 
 is the unique $\mX$-convex viscosity solution of the Dirichlet problem
\begin{equation}
\label{GaussHeisk}
-\det D_{\mathcal X}^2 u + k(x,u)\left(1+|D_{\mathcal X}u|^2\right)^2  = 0 \;\;\mbox{ in }  B_\mathbb{H}(R)%
,\quad
u=R^4\;  \mbox{ on }\; \partial B_\mathbb{H}(R).%
\end{equation}
 In particular, $\ul u$ is the unique $\mX$-convex function on $B_\mathbb{H}(R)$ with 
 horizontal Gauss curvature  $k$ and boundary value $R^4$.
\end{corollary}
\begin{proof}
By Proposition \ref{explicit1} $w(x)=|x|^4_\mathbb{H}$ is a subsolution of \eqref{GaussHeisk} that attains continuously the boundary data.
Then the conclusion follows from Corollary \ref{cnsl} .
\end{proof}
%



\medskip
\noindent
{\bf {\sc  Acknowledgments.}} 
The authors are grateful to Roberto Monti and Luigi Salce for several useful talks and to Luca Capogna for bringing to their attention the paper \cite{RT}.
Work partially supported by the Italian M.I.U.R. project "Viscosity,
metric, and control theoretic methods for nonlinear partial differential equations''.





\end{document}